\documentclass[reqno, 11pt]{amsart}
\usepackage{a4wide}

\usepackage[T1]{fontenc}
\usepackage[utf8]{inputenc}
\usepackage[english]{babel}
\usepackage{amsmath,amssymb,amsfonts}
\usepackage{booktabs}
\usepackage{caption}
\usepackage{mathtools}

\usepackage{pgfplots}
\pgfplotsset{compat=1.18}
\usepgfplotslibrary{fillbetween}
\usetikzlibrary{arrows}
\usepackage{mathrsfs}

\usepackage{bm,cite,xcolor,float,soul,graphicx}
\usepackage[none]{hyphenat}
\usepackage{csquotes}
\usepackage{enumitem}
\usepackage{url}

\usepackage[colorlinks=true,urlcolor=blue,citecolor=violet,linkcolor=teal,linktocpage,pdfpagelabels,bookmarksnumbered,bookmarksopen]{hyperref}
\usepackage{doi}
\usepackage{orcidlink}

\makeatletter
\def\namedlabel#1#2{\begingroup
    #2%
    \def\@currentlabel{#2}%
    \phantomsection\label{#1}\endgroup
}
\makeatother

\usepackage{appendix}
\usepackage{quoting}
\quotingsetup{font=small}
\usepackage{esint}

\usepackage{cleveref}

\usepackage{aliascnt}
\newcommand{\sharedtheorem}[2]{%
  \newaliascnt{#1}{theorem}%
  \newtheorem{#1}[#1]{#2}%
  \aliascntresetthe{#1}%
  \crefname{#1}{#2}{#2s}%
  \Crefname{#1}{#2}{#2s}}
\theoremstyle{plain}
\newtheorem{theorem}{Theorem}[section]
\crefname{theorem}{Theorem}{Theorems}
\Crefname{theorem}{Theorem}{Theorems}
\sharedtheorem{corollary}{Corollary}
\sharedtheorem{lemma}{Lemma}
\sharedtheorem{proposition}{Proposition}
\theoremstyle{definition}
\sharedtheorem{definition}{Definition}
\sharedtheorem{example}{Example}
\sharedtheorem{remark}{Remark}

\title[Selections of IFSs]{Continuous, locally Lipschitz and Lipschitz selections for intuitionistic fuzzy sets, with applications to expert voting}

\author[G. Failla]{Giuseppe Failla\,\orcidlink{0009-0007-7336-4410}}
\address[G. Failla]{Department of Mathematics and Computer Sciences,
Via Archirafi, 34 90123, 
Università degli Studi di Palermo, 
Palermo, Italy}
\email{giuseppe.failla02@unipa.it}
\thanks{Corresponding author: Giuseppe Failla.}

\author[G. La Rosa]{Gianmarco La Rosa\,\orcidlink{0000-0003-1047-5993}}
\address[G. La Rosa]{Department of Mathematics and Computer Sciences,
Via Archirafi, 34 90123, 
Università degli Studi di Palermo, 
Palermo, Italy}
\email{gianmarco.larosa@unipa.it}

\author[M.E. Tabacchi]{Marco Elio Tabacchi\,\orcidlink{0000-0001-7230-8192}}
\address[M.E. Tabacchi]{Department of Mathematics and Computer Sciences,
Via Archirafi, 34 90123, 
Università degli Studi di Palermo, 
Palermo, Italy\\
and
Istituto Nazionale di Ricerche Demopolis
Via Wagner 13, 90139, Palermo, Italy}
\email{marcoelio.tabacchi@unipa.it}


\begin{document}

    \begin{abstract}
	This paper studies continuous selections in the framework of intuitionistic
	fuzzy sets. To each intuitionistic fuzzy set $A=(\mu,\nu)$ on a topological space
	$X$ we associate the interval-valued map $F_{A}(x)=[\mu(x),1-\nu(x)]$, whose width
	is the hesitation degree. Hence, in this setting, a selection of $F_{A}$ is an ordinary fuzzy membership
	function compatible with the intuitionistic information, and every selection is
	the membership function of Atanassov's operator $D_{\vartheta}$ for a suitable
	pointwise parameter $\vartheta$. When $X$ is normal and $\mu$, $\nu$ are upper
	semi-continuous, a continuous selection exists by the Kat\v{e}tov--Tong insertion
	theorem. We show that upper semi-continuity arises naturally in a voting model, in
	which a panel of experts divides a metric universe into closed regions of
	acceptance and rejection, and for this model we construct an explicit continuous
	selection, the consensus selection, illustrated by a non-compact example in the
	plane. On metric spaces we characterise the intuitionistic fuzzy sets that admit
	an $L$-Lipschitz selection through a non-local version of the condition
	$\mu+\nu\leq 1$, and we give the largest and the smallest such selections in
	closed form. Finally, we prove that the set of continuous selections depends in a
	$1$-Lipschitz way on the intuitionistic fuzzy set, with respect to the
	Hausdorff-type distance of Grzegorzewski. We also characterise the intuitionistic fuzzy sets admitting locally Lipschitz selections.
    \end{abstract}
    \maketitle

	\section{Introduction}
	\label{sec:introduction}

	Intuitionistic Fuzzy Sets (IFSs), introduced by Atanassov as a natural
	extension of Zadeh's classical fuzzy sets, provide a flexible framework for modelling
	uncertainty, imprecision, and incomplete information \cite{Ata2,Ata3,Zad}.
	Each element $x\in X$ is associated with both a membership degree
	$\mu(x)\in[0,1]$ and a non-membership degree $\nu(x)\in[0,1]$, subject to the condition
	\[
		0\leq \mu(x)+\nu(x)\leq 1.
	\]
	The remaining quantity
	\[
		\pi(x)=1-\mu(x)-\nu(x)
	\]
	is usually interpreted as the hesitation degree associated with $x$. This
	additional component makes IFSs suitable for settings in which positive and
	negative assessments do not determine the available information completely, as
	frequently occurs in decision-making, expert judgment aggregation, knowledge representation,
	and fuzzy control.

	Fuzzy methods are widely used in the treatment of expert information. Recent approaches
	include interval-based distances for aggregating expert judgments, 
    see Failla--La Rosa--Tabacchi \cite{FaiLaRTab}
	and fuzzy rule-based procedures for weighting heterogeneous groups of human and
	artificial experts, see, e.g., Castronovo and co-authors \cite{CasFilGiaLaRTab, CasFilLaRSanTab}. 
    In these cases the output of the
	procedure is a single value, and its stability under small changes in the input
	is a relevant requirement. A stability result of this kind for continuous
	selections is proved in \Cref{thm:stability}.

	An IFS $A=(\mu,\nu)$ naturally determines the interval
	\[
		F_{A}(x)=[\mu(x),1-\nu(x)]
	\]
	at each point $x\in X$, this correspondence between intuitionistic and
	interval-valued fuzzy sets is classical, see Atanassov--Gargov \cite{AtaGar}. The
	condition defining an IFS guarantees that this interval is nonempty, while its
	width is precisely the hesitation degree $\pi(x)$. A function $g:X\to[0,1]$ satisfying
	\[
		\mu(x)\leq g(x)\leq 1-\nu(x)
	\]
	for every $x\in X$ may therefore be regarded as a single-valued fuzzy membership
	function compatible with the original intuitionistic information. In the
	language of set-valued analysis, such a function is a selection of the interval-valued
	map $F_{A}$. The purpose of this paper is to study whether these selections
	can be chosen continuously, how they can be constructed explicitly, and how they
	depend on the intuitionistic data, under assumptions on $\mu$, $\nu$ and on the
	universe $X$ that are stated explicitly in each result.

	A simple pointwise choice is given by the midpoint of the admissible interval,
	\begin{equation}
		\label{eq:midpoint-introduction}g_{\mathrm{mid}}(x) = \frac{\mu(x)+1-\nu(x)}{2}
		.
	\end{equation}
    
	The midpoint \eqref{eq:midpoint-introduction} is not new: it is the membership
	function of the fuzzy set $D_{1/2}(A)$, where, for $\alpha\in[0,1]$,
	\[
		D_{\alpha}(A)=\bigl\{\bigl(x,\ \mu(x)+\alpha\,\pi(x),\ \nu(x)+(1-\alpha)\,\pi(x)\bigr)\,|\,x\in X\bigr\}
	\]
	is the operator introduced by Atanassov \cite[Section 4.4]{Ata2}. Since the
	hesitation degree of $D_{\alpha}(A)$ vanishes identically, $D_{\alpha}(A)$ is an
	ordinary fuzzy set, and its membership function $(1-\alpha)\mu+\alpha(1-\nu)$
	takes values in $[\mu,1-\nu]$. More generally, letting the parameter depend on
	the point, every function $g$ with $\mu\leq g\leq 1-\nu$ is the membership
	function of $D_{\vartheta}(A)$, defined pointwise as above, for some
	$\vartheta\colon X\to[0,1]$: it suffices to take
	$\vartheta(x)=(g(x)-\mu(x))/\pi(x)$ where $\pi(x)>0$, and $\vartheta(x)$
	arbitrary where $\pi(x)=0$. The problem studied in this paper is therefore the
	choice of the parameter of $D_{\vartheta}$ for which the resulting fuzzy set is
	continuous. For a constant $\alpha$, the membership function of $D_{\alpha}(A)$
	is continuous if and only if $(1-\alpha)\mu-\alpha\nu$ is continuous, which in
	general fails for every $\alpha$: for $\alpha=1/2$ this is the condition that
	$\mu-\nu$ be continuous. Conversely, a continuous selection may require a
	discontinuous parameter $\vartheta$, see \Cref{subsec:voting-example}.

	When $\mu$ and $\nu$ are continuous, the function $g_{\mathrm{mid}}$ is continuous
	as well. Full continuity, however, may be restrictive in models involving
	threshold effects, noisy observations, or abrupt changes between different
	regimes. A natural topological relaxation consists in assuming that the membership
	and non-membership functions are merely upper semi-continuous. Under this
	assumption, the hesitation degree $\pi$ is lower semi-continuous, but the
	midpoint selection need not be continuous. The discontinuities of $\mu$ and $\nu$
	may sometimes cancel each other, although this cannot be expected in general.

	This observation leads directly to a selection problem. By interpreting $F_{A}$
	as a multifunction from $X$ into the family of closed subsets of $\mathbb{R}$,
	the existence of a continuous compatible membership function can be addressed
	through classical results of multivalued analysis. Indeed, if $\mu$ and $\nu$
	are upper semi-continuous, then $1-\nu$ is lower semi-continuous, so that the
	lower endpoint of $F_{A}$ is upper semi-continuous and its upper endpoint is
	lower semi-continuous. Consequently, the interval-valued map $F_{A}$ is lower
	semi-continuous and has nonempty, closed, and convex values. Michael's continuous
	selection theorem then applies on paracompact Hausdorff spaces, while the Kat\v{e}tov--Tong
	insertion theorem gives the corresponding scalar result on normal spaces. These
	results allow us to obtain continuous selections without requiring the membership
	and non-membership functions themselves to be continuous.

	The existence of a continuous selection is a qualitative statement: it
	guarantees compatibility with the intuitionistic information, but it does not control
	how rapidly the selected membership function may vary. On metric domains one
	may ask for the stronger property of a selection satisfying a prescribed
	Lipschitz bound, which would allow its sensitivity to perturbations of the input
	to be estimated quantitatively. This question is addressed in
	\Cref{sec:lipschitz}, where we characterise the intuitionistic fuzzy sets that
	admit an $L$-Lipschitz selection by a non-local version of the inequality
	$\mu+\nu\leq 1$, and we identify the largest and the smallest such selections
	explicitly. Moreover, we characterise, by a non-uniform version of the same inequality,
	the intuitionistic fuzzy sets that admit a locally Lipschitz selection. In the
	same section we show that the family of continuous selections
	depends in a $1$-Lipschitz way on the intuitionistic fuzzy set, when perturbations
	are measured by the Hausdorff-type distance of Grzegorzewski \cite{Grz}.

	Continuous and measurable selections have already appeared in several areas
	related to fuzzy and set-valued analysis. In the fuzzy setting, the seminal
	works of Puri and Ralescu \cite{PurRad1,PurRad2} study measurable selections
	with respect to suitable $\sigma$-algebras in the theory of fuzzy random
	variables. Their measure-theoretic framework is distinct from the topological and
	metric selection problems considered here.

	The interval $F_{A}(x)$ also has a natural geometric connection with the Footprint
	of Uncertainty of an Interval Type-2 Fuzzy Set, see Mendel--John \cite{MenJoh}. In this
	interpretation, $\mu(x)$ and $1-\nu(x)$ play the roles of the lower and upper membership
	functions, respectively, while the hesitation degree $\pi(x)$ gives the width of
	the corresponding footprint. A selection of $F_{A}$ determines an embedded
	type-1 membership function within these bounds. This connection is structural,
	since the selection of a single embedded membership function is conceptually
	different from the type-reduction procedures used in the Interval Type-2
	setting.

	The rest of the paper is organized as follows. \Cref{sec:main_results}
	collects the notions from intuitionistic fuzzy set theory and set-valued analysis
	used throughout, together with the selection and insertion theorems on which our
	arguments rely. \Cref{sec:cont_def} states the basic existence result for upper
	semi-continuous $\mu$ and $\nu$ on normal spaces, which is a direct consequence of
	the Kat\v{e}tov--Tong theorem, discusses the non-uniqueness of selections and the
	degenerate case of vanishing hesitation, and illustrates the construction on
	trapezoidal intuitionistic fuzzy numbers. \Cref{sec:voting} introduces
	intuitionistic fuzzy sets generated by panels of experts on a metric universe, in
	which upper semi-continuity is a consequence of the model rather than an
	assumption; for these sets we construct an explicit continuous selection, the
	consensus selection, and we work out a non-compact example in the plane.
	\Cref{sec:lipschitz} contains the quantitative results: a characterisation of the
	intuitionistic fuzzy sets that admit $L$-Lipschitz selections, with the extremal
	selections in closed form, the corresponding characterisation for locally
	Lipschitz selections, and a stability theorem with respect to the Hausdorff
	distance between intuitionistic fuzzy sets. \Cref{sec:conclusion} discusses
	applications, the relation with interval type-2 fuzzy sets, the reduction of
	other generalizations of intuitionistic fuzzy sets to the present setting, and a
	criterion for choosing a distinguished selection. \Cref{sec:final} contains some
	concluding remarks and open problems.

	\section{Preliminaries}
	\label{sec:main_results}

	Let us recall the definition of an intuitionistic fuzzy set (briefly, IFS), see
	Atanassov \cite{Ata1, Ata2, Ata3}.
	\begin{definition}
		\label{def:IFS} Let $X$ be a nonempty (crisp) set. An \textit{intuitionistic fuzzy
		set in $X$} is an object $A$ given by
		\begin{equation*}
			A = \left\{ \left( x,\, \mu(x),\, \nu(x) \right) \,|\, x\in X \right\},
		\end{equation*}
		where $\mu\colon X \to [0,1]$ and $\nu \colon X \to [0,1]$ satisfy, for every
		$x\in X$, the following inequality
		\begin{equation*}
			0\le \mu(x)+\nu(x)\le 1.
		\end{equation*}
		Here, the functions $\mu$ and $\nu$ represent, respectively, the degree of membership
		and non-membership of an element in $A$. The difference
		$\pi(x) = 1 -\mu(x)-\nu(x)$ represents our lack of knowledge regarding the
		membership and non-membership functions, i.e., the hesitation.
	\end{definition}
	Notice that the geometrical constraint on membership and non-membership ensures
	that the IFS lies in the simplex defined by $0\leq \mu(x)+\nu(x)\leq 1$ for
	all $x\in X$.

    Next, we recall some useful tools form set-valued analysis.
	For a comprehensive introduction on these topics, we refer to Klein--Thompson
	\cite{KleTho}.
	\begin{definition}
		Let $X,Y$ be two sets, we say that the map $F:X\to 2^{Y}$ is a \textit{set-valued
		map} or \textit{multifunction}, i.e.\ $F$ associates with each element of $X$ a subset
		of $Y$.
	\end{definition}
    Notice that in the previous definition no topology is needed.
    In \Cref{def:IFS} the universe $X$ is an
	arbitrary set. Since this paper is concerned with continuity, from now on $X$ is
	a topological space and $\mu$, $\nu$ are real functions on $X$. Further
	assumptions on $X$ (normality, paracompactness, metrizability, compactness) are
	not made globally: each of them is stated explicitly in the results in which it
	is used. In particular, compactness of $X$ is assumed only in
	\Cref{subsec:criterion}.

	\begin{definition}
		Let $X$, $Y$ be topological spaces. Consider a set-valued map $F:X\to 2^{Y}$.
		We say that
		\begin{itemize}
			\item[\textnormal{(i)}] $F$ is \textit{upper semi-continuous} (briefly \textit{u.s.c.}), if for
				every closed subset $C\subset Y$, the set
				\[
					F^{-}(C)=\{x\in X:F(x)\cap C\neq\emptyset\}
				\]
				is closed in $X$. Alternatively, $F$ is u.s.c., if for every open subset
				$O\subset Y$, the set
				\[
					F^{+}(O)=\{x\in X:F(x)\subset O\}
				\]
				is open in $X$.

			\item[\textnormal{(ii)}] $F$ is \textit{lower semi-continuous} (briefly \textit{l.s.c.}), if
				for every closed subset $C\subset Y$, the set
				\[
					F^{+}(C)=\{x\in X:F(x)\subset C\}
				\]
				is closed in $X$. Alternatively, $F$ is l.s.c.\ if for every open subset
				$O\subset Y$, the set
				\[
					F^{-}(O)=\{x\in X:F(x)\cap O\neq\emptyset\}
				\]
				is open in $X$.

			\item[\textnormal{(iii)}] We say that $F$ is \textit{continuous} if it is u.s.c.
				and l.s.c.
		\end{itemize}
	\end{definition}

	The terminology of the above definition should not be confused with the
	semi-continuity of real functions, recalled below. 
	\begin{definition}
		Let $X$ be a topological space. A function
		$f:X\to \mathbb{R}\cup\{\pm\infty\}$ is called
		\begin{itemize}
			\item[\textnormal{(i)}] \textit{upper semi-continuous at} $x_{0}\in X$, if
				for all $\varepsilon>0$ there exists $U$ neighbourhood of $x_{0}$ such
				that
				\[
					f(x_{0})>f(x)-\varepsilon\quad\textnormal{ for all }x\in U;
				\]

			\item[\textnormal{(ii)}] \textit{lower semi-continuous at} $x_{0}\in X$, if
				for all $\varepsilon>0$ there exists $U$ neighbourhood of $x_{0}$ such
				that
				\[
					f(x_{0})<f(x)+\varepsilon\quad\textnormal{ for all }x\in U;
				\]

			\item[\textnormal{(iii)}] \textit{upper} (resp.\ \textit{lower}) \textit{semi-continuous on} $X$
				if it is upper (resp.\ lower) semi-continuous at every point of $X$.
		\end{itemize}
		A function is continuous if and only if it is both upper and lower
		semi-continuous.
	\end{definition}

    The two notions of semi-continuity are related by the 
    following:
    let $F:X\to 2^{\mathbb{R}}$ be a set-valued function and let
	$a,b:X\to\mathbb{R}$ be two functions such that $a(x)\leq b(x)$ for all
	$x\in X$. We call $F$ an interval-valued map if it is defined as
	\[
		F(x)=[a(x),b(x)]\,\,\textnormal{ for all }x\in X.
	\]
    The following is a 
    classical result for interval-valued maps, see Hu--Papageorgiou
	\cite[Example 2.8]{HuPap}.
	\begin{theorem}
		\label{Thm:sciInterval} Let $X\neq \emptyset$ be a topological space and let
		$F:X\to 2^{\mathbb{R}}$ be an interval set-valued map defined as $F(x)=[a(x),
		b(x)]$ with $a,b:X\to \mathbb{R}$ functions such that $a(x)\leq b(x)$ for all
		$x\in X$. Then $F$ is l.s.c. in the sense of set-valued maps if and only if $a$
		is u.s.c. and $b$ is l.s.c. in the sense of functions.
	\end{theorem}

    	We recall the following selection theorem, due to Michael
	\cite[Theorem 3.2]{Mic}.
	\begin{theorem}
		\label{Thm:MichaelSelection} Let $X$ be a paracompact Hausdorff topological space and $Y$
		be a Banach space. Let $F:X\to 2^{Y}$ be a set-valued map l.s.c. and with nonempty,
		closed and convex values. Then $F$ possesses a continuous selection, i.e.\ a continuous
		function $g:X\to Y$, such that $g(x)\in F(x)$ for all $x\in X$.
	\end{theorem}

	We also recall the insertion theorem of Kat\v{e}tov and Tong, see Kat\v{e}tov
	\cite{Kat} and Tong \cite{Tong}; for an extension to functions with values in
	complete lattices see Kubiak--Mardones-P\'erez \cite{KubMar}.
	\begin{theorem}
		\label{thm:KatTon} A topological space $X$ is normal ($T_{4}$) if and only if
		the following holds: for every pair of functions $a,b\colon X\to\mathbb{R}$ such
		that $a$ is upper semi-continuous, $b$ is lower semi-continuous and
		$a(x)\leq b(x)$ for all $x\in X$, there exists a continuous function
		$g\colon X\to\mathbb{R}$ such that
		\[
			a(x) \leq g(x) \leq b(x) \quad \forall x \in X.
		\]
	\end{theorem}

    \begin{remark}
        \label{rem:metricspace-construction} By the proof of Kat\v{e}tov--Tong Theorem,
		if we restrict our contex on metric spaces $X$, we may obtain a formulation for
		our continuous selection $g$. Let $h_{k}:X\to\mathbb{R}$ be continuous
		functions that come from the Urysohn Lemma. Thus, we may write
		\begin{equation}
			\label{Eq:explicit-g}g(x)=\sum_{k=1}^{\infty}2^{-k}\,h_{k}(x),
		\end{equation}
		and
		\[
			g(x)\in[\mu(x),1-\nu(x)]
		\]
		for every $x\in X$. Hence, $g$ is a continuous selection of $F$. We refer the
		interested reader to Engelking \cite[1.7.15., p.62-63]{Eng}, for a
		comprehensive analysis of this case and a detailed proof.
        However, this construction does not provide, in general, a computable closed formula for $g$ in terms of $a$ and $b$, even
		on metric spaces. For the intuitionistic fuzzy sets considered in
		\Cref{sec:voting} an explicit selection is given in \Cref{prop:voting-urysohn},
		and on metric spaces explicit Lipschitz selections are given in
		\Cref{thm:lipschitz}.
    \end{remark}

	\section{Continuous Selections for Intuitionistic Fuzzy Sets}
	\label{sec:cont_def} Inspired by Atanassov \cite[Chapter 2]{Ata2}, let us consider
	the following set-valued map $F:X\to 2^{\mathbb{R}}$ defined as
	\[
		F(x)=[\mu(x),1-\nu(x)]
	\]
	being $\mu,\nu$ as in \Cref{def:IFS}. The geometrical constraint
	$\mu(x)+\nu(x)\leq 1$ ensures that the map $F$ is well defined, with nonempty,
	closed, and convex values.

	The following existence result is a direct consequence of the classical theorems
	recalled in \Cref{sec:main_results}. We state it explicitly, in the form in which
	it will be used in the rest of the paper, and we do not claim originality for it.

	\begin{proposition}
		\label{prop:existence}
		Let $X$ be a normal ($T_{4}$) topological space and let $A=(\mu,\nu)$ be an
		intuitionistic fuzzy set in $X$ such that $\mu$ and $\nu$ are upper
		semi-continuous. Then there exists a continuous function $g\colon X\to[0,1]$
		such that
		\[
			\mu(x)\leq g(x)\leq 1-\nu(x)\qquad\text{for all }x\in X.
		\]
		In particular, the conclusion holds if $X$ is metrizable, or paracompact
		Hausdorff.
	\end{proposition}

	\begin{proof}
		Since $\nu$ is u.s.c., $1-\nu$ is l.s.c., and $\mu\leq 1-\nu$ because $A$ is an
		intuitionistic fuzzy set. By \Cref{thm:KatTon} there exists a continuous
		$g\colon X\to\mathbb{R}$ with $\mu\leq g\leq 1-\nu$. Since $0\leq\mu$ and
		$1-\nu\leq 1$, the function $g$ takes values in $[0,1]$. Every metrizable space
		and every paracompact Hausdorff space is normal, see Engelking \cite{Eng}.
	\end{proof}

	\begin{remark}
		\label{rem:michael-route}
		If $X$ is paracompact Hausdorff, the same conclusion can be obtained from
		Michael's theorem: by \Cref{Thm:sciInterval}, $F$ is l.s.c.\ as a set-valued map,
		and it has nonempty, compact and convex values, so that
		\Cref{Thm:MichaelSelection} applies. Since paracompact Hausdorff spaces are
		normal, this route yields nothing more than \Cref{prop:existence} for
		real-valued selections; it is, however, the one that extends to set-valued maps
		with values in higher-dimensional spaces, see \Cref{subsec:other-sets}.
	\end{remark}

		Clearly, one may wonder whether it makes sense to assume that $\mu$ and
		$\nu$ are u.s.c. The following considerations can give the answer, and it is
		at least twofold. First, a \enquote{principle of conservativeness and caution}
		could be observed. Consider a point approaching a critical value $x_{0}$. The
		membership degree $\mu$ at $x_{0}$ cannot be smaller than the function behaviour
		of the surrounding values. Consequently, the model never underestimates the certainty
		associated with an exact observation. Similarly, the non-membership degree $\nu$
		at $x_{0}$ cannot be smaller than the value of the surrounding non-membership
		values. Therefore, the model does not arbitrarily reduce rejection or uncertainty
		at boundary points.
		On the other hand, the hesitation degree $\pi=1-\mu-\nu$ is lower
		semi-continuous, since $-\mu$ and $-\nu$ are lower semi-continuous. Thus, for every
		$\varepsilon>0$ there exists a neighbourhood $U$ of $x_{0}$ such that
		$\pi(x_{0})\leq\pi(x)+\varepsilon$ for all $x\in U$, that is,
		$\pi(x_{0})\leq\liminf\limits_{x\to x_{0}}\pi(x)$. In words, the hesitation at a point
		cannot exceed the hesitation observed at nearby points: an exact observation is
		at least as informative as the observations in its vicinity, and hesitation
		cannot appear abruptly at an isolated point. Lower semi-continuity does not
		prevent the opposite phenomenon: the hesitation may drop at a point where it is
		positive nearby. This happens, for instance, in the voting model of
		\Cref{sec:voting}, on the boundary of the region where an expert abstains; in the
		example of \Cref{subsec:voting-example}, $\pi$ vanishes on two circles and equals
		$\frac{1}{3}$ in the surrounding annuli.

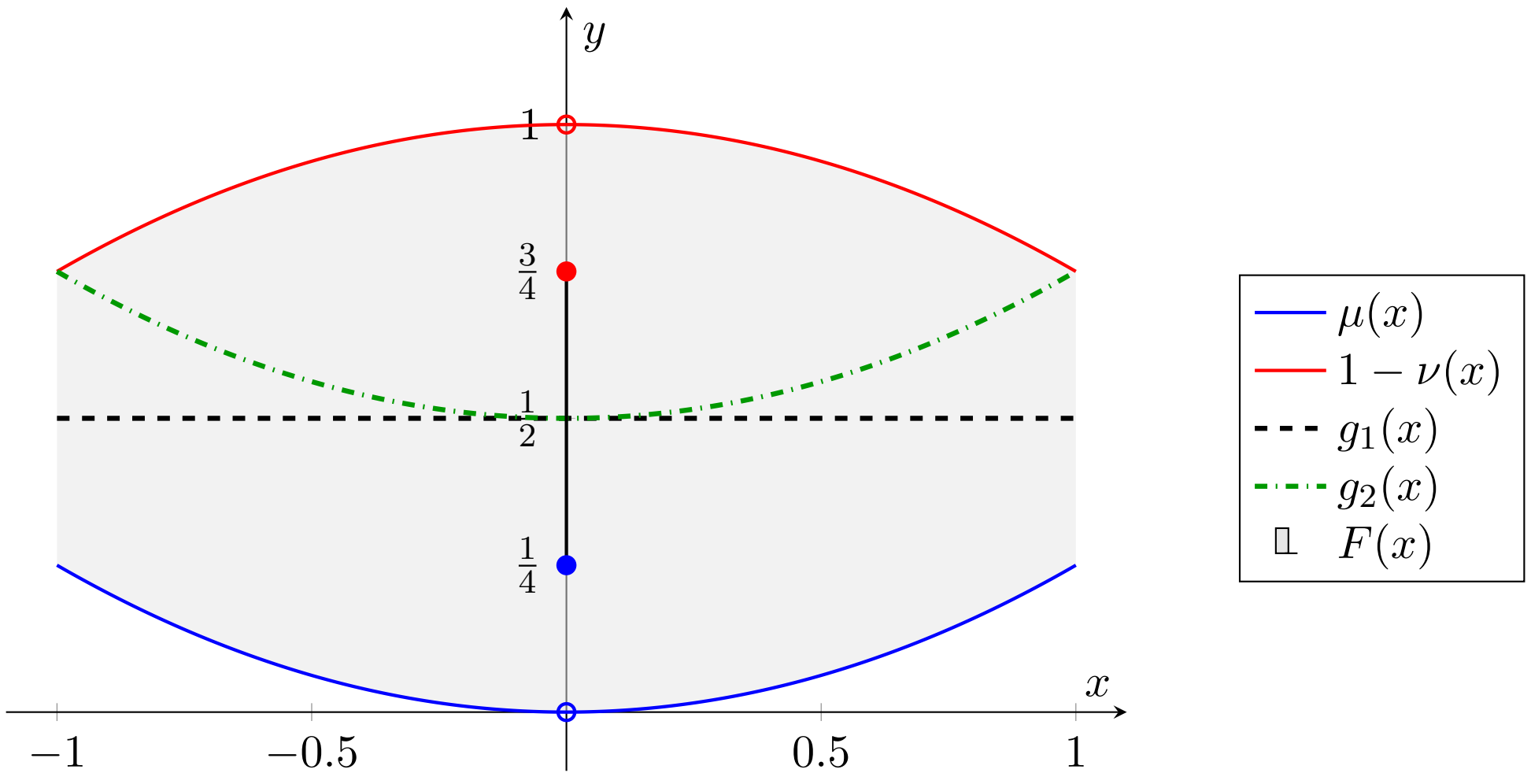
\begin{figure}[ht]
		\centering
		\begin{tikzpicture}[scale=1]
			\begin{axis}[
				width=11cm,
				height=8cm,
				axis lines=middle,
				xmin=-1.1,
				xmax=1.1,
				ymin=-0.1,
				ymax=1.2,
				xlabel={$x$},
				ylabel={$y$},
				xtick={-1, -0.5, 0, 0.5, 1},
				ytick={0, 0.25, 0.5, 0.75, 1},
				yticklabels={$0$, $\frac{1}{4}$, $\frac{1}{2}$, $\frac{3}{4}$, $1$},
				legend style={at={(1.1,0.65)}, anchor=north west},
				legend cell align=left,
				samples=200
			]
				\addplot[blue, thick, domain=-1:1] {x^2/4}; \addlegendentry{$\mu(x)$}
				\addplot[red, thick, domain=-1:1] {1 - x^2/4}; \addlegendentry{$1-\nu(x)$}
				\addplot[black, dashed, very thick, domain=-1:1] {0.5}; \addlegendentry{$g_{1}(x)$}
				\addplot[green!60!black, dashdotted, very thick, domain=-1:1]
				{0.5 + x^2/4}; \addlegendentry{$g_{2}(x)$}
				\addlegendimage{ybar, xbar legend, fill=gray!20, draw=none}
				\addlegendentry{$F(x)$}
				\addplot[name path=lower, draw=none, domain=-1:1, forget plot] {x^2/4};
				\addplot[name path=upper, draw=none, domain=-1:1, forget plot] {1 - x^2/4};
				\addplot[gray!20, opacity=0.5, forget plot] fill between[of=lower and
				upper];
				\addplot[mark=o, blue, thick, fill=white, forget plot] coordinates {(0,0)};
				\addplot[mark=o, red, thick, fill=white, forget plot] coordinates {(0,1)};
				\addplot[mark=*, blue, thick, forget plot] coordinates {(0,0.25)}; \addplot[mark=*,
				red, thick, forget plot] coordinates {(0,0.75)};
				\draw[black, thick] (axis cs:0,0.25) -- (axis cs:0,0.75);
			\end{axis}
		\end{tikzpicture}
		\caption{The set-valued map $F$ of \Cref{Example1}, its endpoints and the two
		continuous selections $g_{1}$ and $g_{2}$.}
		\label{fig:example1}
	\end{figure}

	\begin{example}
		\label{Example1} Define $\mu, \nu: [-1, 1] \to [0,1]$ as follows:
		\begin{align}
			\mu(x) & = \begin{cases}\frac{x^{2}}{4} & \text{if } x \neq 0, \\ \frac{1}{4} & \text{if } x = 0,\end{cases} \qquad \nu(x) = \begin{cases}\frac{x^2}{4}&\text{if }x \neq 0, \\ \frac{1}{4}&\text{if }x = 0.\end{cases}
		\end{align}
		Let us explicitly evaluate the multifunction $F(x) = [\mu(x), 1 - \nu(x)]$:
		\begin{itemize}
			\item For $x \neq 0$:
				\begin{equation}
					F(x) = \left[ \frac{x^{2}}{4}, 1 - \frac{x^{2}}{4}\right].
				\end{equation}

			\item For $x = 0$:
				\begin{equation}
					F(0) = \left[ \frac{1}{4}, 1 - \frac{1}{4}\right] = \left[ \frac{1}{4},
					\frac{3}{4}\right].
				\end{equation}
		\end{itemize}
		Geometrically, for $x \neq 0$, $F(x)$ forms a ``tube'' that widens as $x \to
		0$, approaching the interval $[0, 1]$. At precisely $x = 0$, the interval
		abruptly constricts to $[1/4, 3/4]$. This sudden constriction, which remains
		strictly within the limiting boundary of the tube, characterises $F(x)$ as lower
		semi-continuous.

		\noindent
		\textbf{Option A: A Constant Selection} \\ We can choose a constant value
		that fits within the narrowest ``bottleneck'' of the set-valued map. Let:
		\begin{equation}
			g_{1}(x) = \frac{1}{2}.
		\end{equation}
		We verify that $g_{1}(x) \in F(x)$ for all $x$:
		\begin{itemize}
			\item At $x = 0$: $1/2 \in [1/4, 3/4]$, which is true.

			\item For $x \neq 0$: Since $x \in [-1, 1]$, we have $x^{2}/4 \leq 1/4$.
				Thus, the lower bound is at most $1/4$, and the upper bound is at least
				$3/4$. The interval $\left[\frac{x^{2}}{4}, 1 - \frac{x^{2}}{4}\right]$ clearly
				contains $1/2$.
		\end{itemize}

		\vspace{0.5cm}
		\noindent
		\textbf{Option B: A Variable Selection} \\ We can construct a non-constant selection
		$g(x)$ that dynamically follows the boundaries of the interval:
		\begin{equation}
			g_{2}(x) = \frac{1}{2}+ \frac{x^{2}}{4}.
		\end{equation}
		The map $F$ and the two selections $g_{1}$, $g_{2}$ are represented in \Cref{fig:example1}.

	\end{example}

	\begin{remark}
		\label{rem:example1-midpoint}
		In \Cref{Example1} the discontinuities of $\mu$ and $\nu$ at $x=0$ compensate
		each other: the midpoint \eqref{eq:midpoint-introduction} is
		$g_{\mathrm{mid}}\equiv\frac{1}{2}$, which is continuous, and a constant selection
		exists because $\sup\limits_{X}\mu+\sup\limits_{X}\nu\leq 1$ (compare with the discussion of
		$L^{*}(A)$ after \Cref{thm:lipschitz}). The example therefore illustrates the
		non-uniqueness of continuous selections, but not the need for
		\Cref{prop:existence}. An example in which the midpoint is discontinuous and no
		constant selection exists is given in \Cref{subsec:voting-example}.
	\end{remark}

	Notice that the above example exhibits two distinct continuous selections of
	the same map; in general, the selection $g$ is not unique. This is not a defect
	of the construction. At each point the set of admissible values has width
	$\pi(x)$, so a single selection can be singled out only by a criterion external
	to the topological data. When the hesitation index vanishes, i.e.,
		$\pi(x) = 0$, the interval of valid selections geometrically collapses,
		yielding a unique continuous defuzzifier $g(x) = \mu(x)$. By definition, the
		hesitation margin in an Intuitionistic Fuzzy Set (IFS) is formulated as
		$\pi(x) = 1 - \mu(x) - \nu(x)$. Setting this parameter to zero implies that $1
		- \nu(x) = \mu(x)$, causing the set-valued mapping
		$F(x) = [\mu(x), 1-\nu(x)]$ to degenerate into the singleton
		$[\mu(x), \mu(x)]$ for all $x\in X$. Consequently, the selection mapping is
		constrained to a single degree of freedom, trivially resolving to $g(x) = \mu
		(x)$.

		Note also that in this case, under the assumption that $\mu$ and $\nu$ are
		u.s.c., $\mu$ is simultaneously u.s.c. and l.s.c., hence
		continuous: the regularity of the membership function is forced by the vanishing
		of the hesitation. Semantically, the condition $\pi\equiv 0$ removes the intuitionistic
		component altogether and the model reduces to a Zadeh fuzzy set \cite{Zad},
		in which non-membership is the complement of membership. The pointwise version
		of the statement also holds: on the set $\{x\in X:\pi(x)=0\}$ every
		selection coincides with $\mu$. The multiplicity of admissible selections is therefore the analytic counterpart
	of the hesitation degree, and choosing among them requires additional
	modelling assumptions. We return to this question in \Cref{subsec:criterion}.

	\begin{example}
		For the sake of clarity, we add a geometrical visualization in the special
		case of trapezoidal intuitionistic fuzzy numbers (TrIFNs), i.e.\ when membership and non-membership
		possess the form
		\begin{equation}
			\label{eq:TrIFN}\mu(x) =
			\begin{cases}
				\dfrac{x-a_1}{a_2 - a_1} & \textnormal{if}\quad a_{1} \le x \le a_{2} \\[0.6em]
				1                        & \textnormal{if}\quad a_{2} \le x \le a_{3} \\[0.6em]
				\dfrac{a_4-x}{a_4 - a_3} & \textnormal{if}\quad a_{3} \le x \le a_{4} \\[0.6em]
				0                        & \textnormal{otherwise}
			\end{cases}
			\,,\quad \nu(x)=
			\begin{cases}
				\dfrac{b_2-x}{b_2 - b_1} & \textnormal{if}\quad b_{1} \le x \le b_{2} \\[0.6em]
				0                        & \textnormal{if}\quad b_{2} \le x \le b_{3} \\[0.6em]
				\dfrac{x-b_3}{b_4 - b_3} & \textnormal{if}\quad b_{3} \le x \le b_{4} \\[0.6em]
				1                        & \textnormal{otherwise}.
			\end{cases}
		\end{equation}
		with $a_{1}\le a_{2}\le a_{3}\le a_{4}$, $b_{1}\le b_{2}\le b_{3}\le b_{4}$ and
		$b_{1} \le a_{1},\,\, b_{2} \le a_{2} \le a_{3} \le b_{3},\,\, a_{4} \le b_{4}$;
		these conditions ensure that $\mu+\nu\leq 1$ on $\mathbb{R}$.

		Note that $\mu$ and $\nu$ are continuous, so that the midpoint
		$g_{\mathrm{mid}}$ is itself a continuous selection and \Cref{prop:existence} is
		not needed here: the example only serves to visualize the set-valued map $F$. We
		also point out that in the literature on intuitionistic fuzzy numbers the
		functions $\mu$ and $1-\nu$ are frequently required to be upper
		semi-continuous, that is, $\nu$ is required to be lower semi-continuous. Under
		that convention both endpoints of $F$ are u.s.c., and the existence of a
		continuous selection is not guaranteed by \Cref{prop:existence}; for TrIFNs the
		two conventions agree, since $\mu$ and $\nu$ are continuous.

		Let us consider a membership and non-membership defined as follows,
		\[
			\mu(x)=
			\begin{cases}
				0,              & x\le 2,   \\[2mm]
				\dfrac{x-2}{2}, & 2<x\le 4, \\[2mm]
				1,              & 4<x\le 6, \\[2mm]
				\dfrac{8-x}{2}, & 6<x\le 8, \\[2mm]
				0,              & x>8,
			\end{cases}
			\quad\textnormal{ and }\quad \nu(x)=
			\begin{cases}
				1,              & x\le 1,   \\[2mm]
				\dfrac{3-x}{2}, & 1<x\le 3, \\[2mm]
				0,              & 3<x\le 7, \\[2mm]
				\dfrac{x-7}{2}, & 7<x\le 9, \\[2mm]
				1,              & x>9.
			\end{cases}
		\]
		It is straightforward to verify that $\mu(x)+\nu(x)\le 1$ for every
		$x\in\mathbb{R}$. The functions $\mu$, $\nu$, the set-valued map $F$
		and the midpoint selection are represented in \Cref{fig:trifn}.
		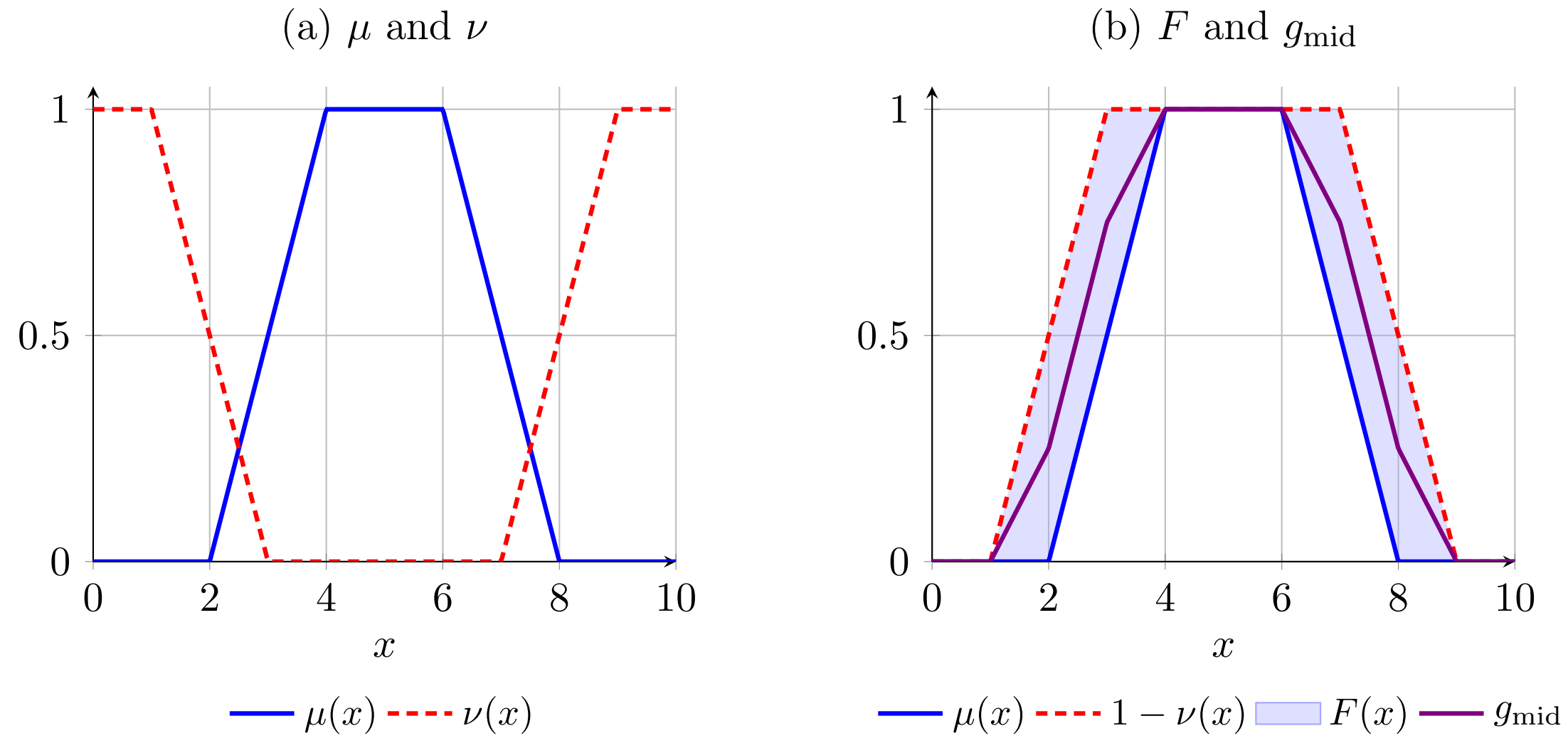
\begin{figure}[ht]
			\centering
			\begin{tikzpicture}
				\begin{axis}[
					width=7cm,
					height=6cm,
					xmin=0,
					xmax=10,
					ymin=0,
					ymax=1.05,
					xlabel={$x$},
					ylabel={},
					title={(a) $\mu$ and $\nu$},
					axis lines=left,
					samples=400,
					domain=0:10,
					grid=major,
					legend style={ draw=none, fill=none, font=\small, at={(0.5,-0.25)}, anchor=north, legend columns=2 }
				]
					\addplot[ blue, very thick] expression{ (x<=2)*0 +(x>2)*(x<=4)*((x-2)/2) +(x>4)*(x<=6)*1 +(x>6)*(x<=8)*((8-x)/2) +(x>8)*0 };

					\addplot[ red, very thick, densely dashed] expression{ (x<=1)*1 +(x>1)*(x<=3)*((3-x)/2) +(x>3)*(x<=7)*0 +(x>7)*(x<=9)*((x-7)/2) +(x>9)*1 };

					\legend{$\mu(x)$,$\nu(x)$}
				\end{axis}

				\begin{axis}[
					at={(7.8cm,0)},
					anchor=south west,
					width=7cm,
					height=6cm,
					xmin=0,
					xmax=10,
					ymin=0,
					ymax=1.05,
					xlabel={$x$},
					title={(b) $F$ and $g_{\mathrm{mid}}$},
					axis lines=left,
					samples=400,
					domain=0:10,
					grid=major,
					legend style={ draw=none, fill=none, font=\small, at={(0.5,-0.25)}, anchor=north, legend columns=4 }
				]
					\addplot[ name path=A, blue, very thick] expression{ (x<=2)*0 +(x>2)*(x<=4)*((x-2)/2) +(x>4)*(x<=6)*1 +(x>6)*(x<=8)*((8-x)/2) +(x>8)*0 };

					\addplot[ name path=B, red, very thick, densely dashed] expression{ 1- ( (x<=1)*1 +(x>1)*(x<=3)*((3-x)/2) +(x>3)*(x<=7)*0 +(x>7)*(x<=9)*((x-7)/2) +(x>9)*1 ) };

					\addplot[ blue!35, fill opacity=0.35] fill between[ of=A and B];

					\addplot[ very thick, violet] expression{ ( ( (x<=2)*0 +(x>2)*(x<=4)*((x-2)/2) +(x>4)*(x<=6)*1 +(x>6)*(x<=8)*((8-x)/2) +(x>8)*0 ) + ( 1- ( (x<=1)*1 +(x>1)*(x<=3)*((3-x)/2) +(x>3)*(x<=7)*0 +(x>7)*(x<=9)*((x-7)/2) +(x>9)*1 ) ) )/2 };

					\legend{$\mu(x)$,$1-\nu(x)$,$F(x)$,$g_{\mathrm{mid}}$}
				\end{axis}
			\end{tikzpicture}
			\caption{The TrIFN of \Cref{eq:TrIFN} with $a=(2,4,6,8)$ and
			$b=(1,3,7,9)$: (a) membership and non-membership functions; (b) the
			set-valued map $F$ and the midpoint selection $g_{\mathrm{mid}}$.}
			\label{fig:trifn}
		\end{figure}
	\end{example}

	\begin{remark}
		We emphasize that, for TrIFNs it is always possible to construct an
		explicit continuous selection $g\in F(x)$ of the following form
		\begin{equation*}
			g(x)=
			\begin{cases}
				0                             & \textnormal{ for }x\leq b_{1},     \\
				\frac{x-b_{1}}{a_{2}-b_{1}}   & \textnormal{ for }b_{1}\leq x\leq a_{2}, \\
				1                             & \textnormal{ for }a_{2}\leq x\leq a_{3}, \\
				1-\frac{x-a_{3}}{b_{4}-a_{3}} & \textnormal{ for }a_{3}\leq x\leq b_{4}, \\
				0                             & \textnormal{ for }x\geq b_{4},
			\end{cases}
		\end{equation*}
		namely, the membership function of the trapezoidal fuzzy number
		$(b_{1},a_{2},a_{3},b_{4})$;
		that is, the function represented in \Cref{FigSelCon}.
		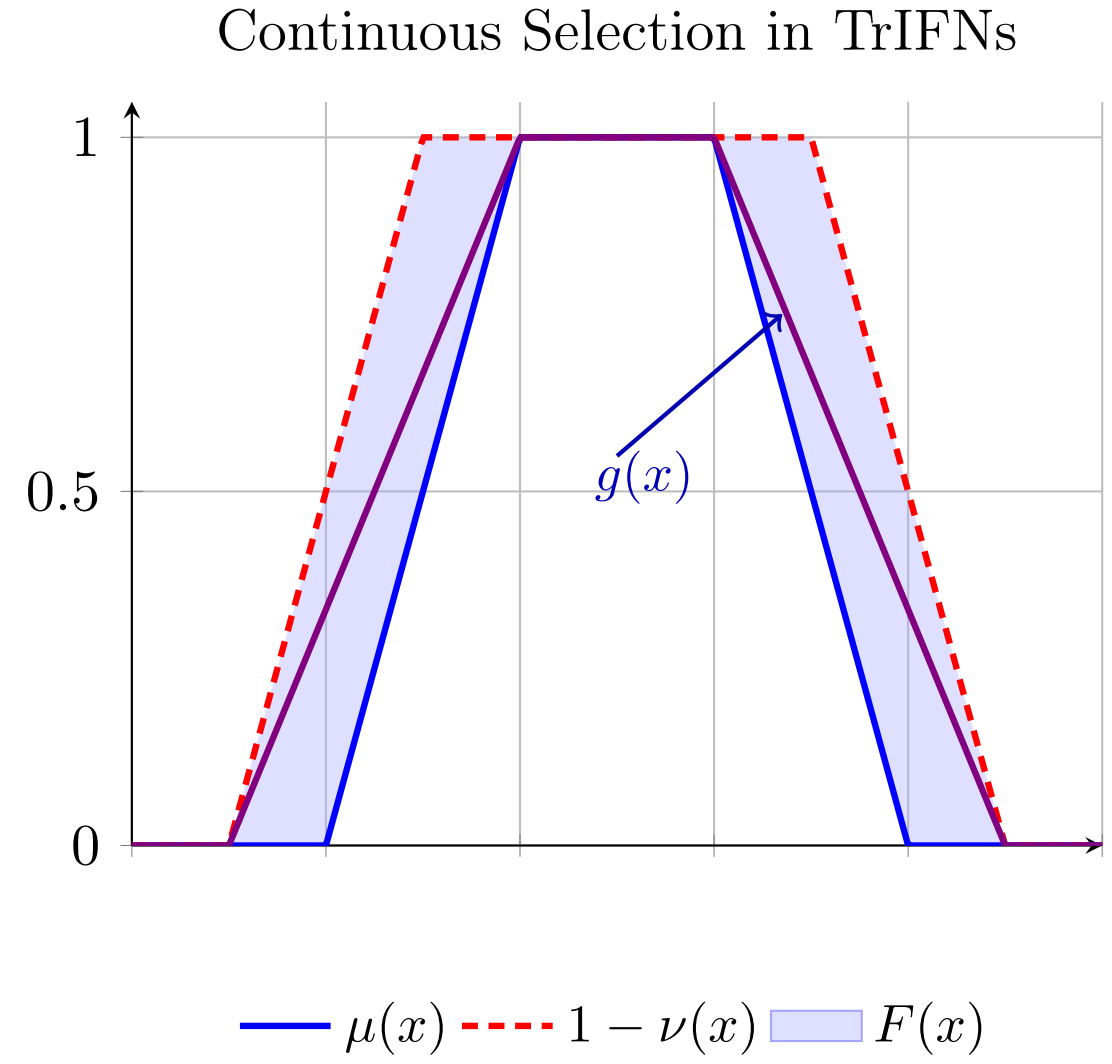
\begin{figure}[h]
			\centering
			\begin{tikzpicture}
				\begin{axis}[
					at={(9.5cm,0)},
					anchor=south west,
					width=8cm,
					height=6.5cm,
					xmin=0,
					xmax=10,
					ymin=0,
					ymax=1.05,
					xlabel={\empty},
					xticklabels={},
					title={Continuous Selection in TrIFNs},
					axis lines=left,
					samples=400,
					domain=0:10,
					grid=major,
					legend style={ draw=none, fill=none, font=\small, at={(0.5,-0.18)}, anchor=north, legend columns=3 }
				]
					\addplot[ name path=A, blue, very thick] expression{ (x<=2)*0 +(x>2)*(x<=4)*((x-2)/2) +(x>4)*(x<=6)*1 +(x>6)*(x<=8)*((8-x)/2) +(x>8)*0 };

					\addplot[ name path=B, red, very thick, densely dashed] expression{ 1- ( (x<=1)*1 +(x>1)*(x<=3)*((3-x)/2) +(x>3)*(x<=7)*0 +(x>7)*(x<=9)*((x-7)/2) +(x>9)*1 ) };

					\addplot[ blue!35, fill opacity=0.35] fill between[ of=A and B];

					\addplot[ very thick, violet] expression{ (x<1)*0 +(x>=1)*(x<=4)*((x-1)/3) +(x>4)*(x<6)*1 +(x>=6)*(x<=9)*((9-x)/3) +(x>9)*0 };
					\node[blue!70!black,font=\small] at (axis cs:5.28,0.52) {$g(x)$};
					\draw[ ->, blue!70!black, thick] (axis cs:5,0.55) -- (axis cs:6.7,0.75);
					\legend{$\mu(x)$,$1-\nu(x)$,$F(x)$}
				\end{axis}
			\end{tikzpicture}
			\caption{Geometric visualization of continuous selection in the TrIFN
			case.}\label{FigSelCon}
		\end{figure}
	\end{remark}

\section{Intuitionistic fuzzy sets generated by expert panels}
\label{sec:voting}

The examples discussed so far were chosen for their simplicity: the universe
was a subset of the real line and the upper semi-continuity of $\mu$
and $\nu$ was imposed directly. In this section we consider a class of
intuitionistic fuzzy sets in which upper semi-continuity is not an assumption
but a consequence of the way the set is constructed. The construction is the
voting interpretation of intuitionistic fuzzy sets, already discussed by
Atanassov \cite[Sections 1.1 and 5]{Ata2}, see also \cite{Ata4} and \cite{Grz}: the membership degree is the proportion of voters in
favour, the non-membership degree is the proportion of voters against, and the
hesitation degree is the proportion of abstentions. The universe is an arbitrary
metric space, which need be neither compact nor a subset of $\mathbb{R}$.

\subsection{The voting model}
\label{subsec:voting-model}

Let $(X,d)$ be a metric space and consider a panel of $n\geq 1$ experts. Each
expert $k\in\{1,\dots,n\}$ partitions the universe into three regions: a set
$A_{k}\subseteq X$ of elements that the expert accepts, a set $R_{k}\subseteq X$ of
elements that the expert rejects, and the remaining set
$X\setminus(A_{k}\cup R_{k})$, on which the expert abstains. We require
\begin{equation}
	\label{eq:voting-hyp}
	A_{k},\,R_{k}\ \text{closed in }X,\qquad A_{k}\cap R_{k}=\emptyset,
	\qquad k=1,\dots,n.
\end{equation}
Closedness expresses the fact that a judgement which holds at every element
of a convergent sequence is retained at the limit point: an expert who accepts
every site at distance strictly less than $r$ from a reference point also accepts
the sites at distance exactly $r$. Disjointness simply means that no expert
accepts and rejects the same element. The panel induces the functions
\begin{equation}
	\label{eq:voting-mu-nu}
	\mu(x)=\frac{1}{n}\sum_{k=1}^{n}\mathbf{1}_{A_k}(x),\qquad
	\nu(x)=\frac{1}{n}\sum_{k=1}^{n}\mathbf{1}_{R_k}(x),\qquad x\in X,
\end{equation}
where $\mathbf{1}_{E}$ denotes the characteristic function of $E\subseteq X$.

\begin{proposition}
	\label{prop:voting-IFS}
	Under \eqref{eq:voting-hyp}, the pair $(\mu,\nu)$ defined in
	\eqref{eq:voting-mu-nu} is an intuitionistic fuzzy set in $X$, the functions
	$\mu$ and $\nu$ are upper semi-continuous, and
	\[
		\pi(x)=\frac{1}{n}\,\#\bigl\{k:\ x\notin A_{k}\cup R_{k}\bigr\}
	\]
	is the proportion of experts abstaining on $x$. Consequently, the set-valued
	map $F(x)=[\mu(x),1-\nu(x)]$ admits a continuous selection.
\end{proposition}

\begin{proof}
	The indicator function of a closed set $E$ is upper semi-continuous, since
	$\{x:\mathbf{1}_{E}(x)\geq t\}$ equals $X$ for $t\leq 0$, $E$ for
	$0<t\leq 1$, and $\emptyset$ for $t>1$, and is therefore closed for every
	$t\in\mathbb{R}$. A finite sum of upper semi-continuous functions with
	nonnegative coefficients is upper semi-continuous, so $\mu$ and $\nu$ are
	upper semi-continuous. Since $A_{k}\cap R_{k}=\emptyset$, we have
	$\mathbf{1}_{A_k}+\mathbf{1}_{R_k}=\mathbf{1}_{A_k\cup R_k}\leq 1$, and
	summing over $k$ gives $\mu+\nu\leq 1$. Thus the formula for $\pi$ follows from
	$1-\mathbf{1}_{A_k\cup R_k}=\mathbf{1}_{X\setminus(A_k\cup R_k)}$. Finally,
	every metric space is normal (see Engelking \cite{Eng}), so the existence of
	a continuous selection follows from \Cref{prop:existence}.
\end{proof}

    Notice that $\pi(x)=0$ exactly when every expert expresses an opinion on $x$. On
    this set every selection is forced to coincide with $\mu$, 
    see the discussion at the end of the previous section. 
    The selections are therefore determined by the panel
    where nobody abstains, and the freedom in choosing them is confined to the region
    where at least one expert abstains. We also observe that the hesitation degree
    of a voting model may vanish at a point while being positive at every nearby
    point, which happens on the boundary of the region where an expert abstains;
    this is compatible with the lower semi-continuity of $\pi$.
    
    \subsection{An explicit selection}
    \label{subsec:voting-explicit}
    
    \Cref{prop:existence} guarantees existence, and
    \Cref{rem:metricspace-construction} explains why the classical proof does not
    provide a closed formula. For the voting model, the
    structure of the panel yields a selection in closed form. For each expert $k$
    define $u_{k}\colon X\to[0,1]$ by the TOPSIS-type function
    \begin{equation}
    	\label{eq:voting-urysohn}
    	u_{k}(x)=\frac{\operatorname{dist}(x,R_{k})}{\operatorname{dist}(x,A_{k})+\operatorname{dist}(x,R_{k})}
    	\qquad\text{if }A_{k}\neq\emptyset\neq R_{k},
    \end{equation}
    and set $u_{k}\equiv 1$ if $R_{k}=\emptyset$, and $u_{k}\equiv 0$ if $R_{k}\neq\emptyset$
    and $A_{k}=\emptyset$. Here
    $\operatorname{dist}(x,E)=\inf\{d(x,e):\ e\in E\}$.
    The function $u_{k}$ is the Urysohn function of the pair $(A_{k},R_{k})$: it
    replaces the three-valued judgement of the $k$-th expert by a continuous
    judgement which equals $1$ on the accepted region, $0$ on the rejected region,
    and interpolates between them on the region of abstention.

\begin{proposition}
	\label{prop:voting-urysohn}
	Under \eqref{eq:voting-hyp}, the function
	\begin{equation}
		\label{eq:voting-consensus}
		g_{P}(x)=\frac{1}{n}\sum_{k=1}^{n}u_{k}(x),\qquad x\in X,
	\end{equation}
	is a continuous selection of $F(x)=[\mu(x),1-\nu(x)]$.
\end{proposition}

\begin{proof}
	Let $k$ be such that $A_{k}\neq\emptyset\neq R_{k}$. The functions
	$x\mapsto\operatorname{dist}(x,A_{k})$ and $x\mapsto\operatorname{dist}(x,R_{k})$
	are continuous, and their sum never vanishes: since $A_{k}$ and $R_{k}$ are
	closed, $\operatorname{dist}(x,A_{k})=0$ if and only if $x\in A_{k}$, and
	similarly for $R_{k}$, while $A_{k}\cap R_{k}=\emptyset$. Hence $u_{k}$ is
	continuous, takes values in $[0,1]$, equals $1$ on $A_{k}$ and equals $0$ on
	$R_{k}$. The same properties hold trivially in the remaining cases. Therefore
	\[
		\mathbf{1}_{A_k}\leq u_{k}\leq 1-\mathbf{1}_{R_k}\qquad\text{on }X,
	\]
	and averaging over $k$ gives $\mu\leq g_{P}\leq 1-\nu$. As a finite average
	of continuous functions, $g_{P}$ is continuous.
\end{proof}

We call $g_{P}$ the \emph{consensus selection} of the panel. It is obtained by
first making the opinion of each expert continuous and then averaging, whereas
the midpoint $g_{\mathrm{mid}}$ averages the discontinuous opinions directly.
In comparison with \Cref{rem:metricspace-construction}, the infinite series of
unspecified continuous functions is replaced by a finite average of $n$
explicit functions, one for each expert.
We also observe that $u_{k}$ has the same form as the relative closeness
coefficient of the TOPSIS method of Hwang and Yoon \cite{HwaYoo}, with the accepted
region $A_{k}$ in the role of the ideal solution and the rejected region $R_{k}$ in the
role of the anti-ideal one: $u_{k}(x)$ is the relative closeness of $x$ to $A_{k}$. The
consensus selection is therefore the average of the TOPSIS scores of the experts,
each computed with respect to the expert's own acceptance and rejection regions.

\begin{remark}
	\label{rem:voting-non-uniform}
	The consensus selection is continuous, but it need not be uniformly
	continuous. In the following example, in fact, no selection is uniformly
	continuous. Let $X=\mathbb{R}^{2}$
	with the Euclidean distance, $n=1$, and
	\[
		A_{1}=\{(x_{1},x_{2}):\ x_{1}>0,\ x_{1}x_{2}\geq 1\},\qquad
		R_{1}=\{(x_{1},x_{2}):\ x_{2}\leq 0\}.
	\]
	Both sets are closed and disjoint, so \Cref{prop:voting-urysohn} applies.
	Since $\pi=0$ on $A_{1}\cup R_{1}$, every selection $g$ equals $1$ on $A_{1}$
	and $0$ on $R_{1}$. The points $(t,1/t)\in A_{1}$ and $(t,0)\in R_{1}$ are at
	distance $1/t$, while $|g(t,1/t)-g(t,0)|=1$ for every $t>0$. Hence no selection
	is uniformly continuous. This phenomenon cannot occur on a compact universe,
	where every continuous function is uniformly continuous, and it shows that
	the non-compact setting is not a mere formal generalization. In particular,
	no selection is Lipschitz, that is, $L^{*}(A)=+\infty$ in the notation of
	\Cref{sec:lipschitz}, although $\mu$ and $\nu$ are u.s.c.; see
	\Cref{rem:independence}.
\end{remark}

\subsection{A worked example: site assessment in the plane}
\label{subsec:voting-example}

We now work out a complete example on a non-compact universe of dimension two.
A panel of $n=3$ experts assesses the suitability of locations
$x\in X=\mathbb{R}^{2}$ for a facility, for instance a service that should be
close to a reference point placed at the origin. Distances are Euclidean and
are measured in kilometres. We write $r=|x|$. The judgements of the experts are
\[
	\begin{array}{c|c|c}
		\text{expert } k & A_{k}\ (\text{accept}) & R_{k}\ (\text{reject}) \\
		\hline
		1 & \{r\leq 2\} & \{r\geq 3\} \\
		2 & \{r\leq 1\} & \{r\geq 2\} \\
		3 & \{r\leq 3\} & \{r\geq 4\}
	\end{array}
\]
All the sets are closed and $A_{k}\cap R_{k}=\emptyset$, so \eqref{eq:voting-hyp}
holds. By \eqref{eq:voting-mu-nu}, $\mu$ and $\nu$ are radial step functions:
\[
	\mu(x)=
	\begin{cases}
		1, & r\leq 1, \\[1mm]
		\frac{2}{3}, & 1<r\leq 2, \\[1mm]
		\frac{1}{3}, & 2<r\leq 3, \\[1mm]
		0, & r>3,
	\end{cases}
	\qquad
	\nu(x)=
	\begin{cases}
		0, & r<2, \\[1mm]
		\frac{1}{3}, & 2\leq r<3, \\[1mm]
		\frac{2}{3}, & 3\leq r<4, \\[1mm]
		1, & r\geq 4,
	\end{cases}
\]
which are upper semi-continuous by \Cref{prop:voting-IFS}, and the hesitation
degree is
\[
	\pi(x)=
	\begin{cases}
		0, & r\leq 1,\ \ r=2,\ \ r=3,\ \ r\geq 4, \\[1mm]
		\frac{1}{3}, & \text{otherwise}.
	\end{cases}
\]
Thus the panel is unanimous on the disc $\{r\leq 1\}$ and outside the disc
$\{r<4\}$, and it is complete, although not unanimous, on the two circles $r=2$
and $r=3$. On each open annulus between consecutive radii, exactly one expert
abstains.

\paragraph{Failure of the elementary choices.}
No constant function is a selection, see again the discussion at the end of \Cref{sec:cont_def}, every
selection equals $1$ on $\{r\leq 1\}$ and $0$ on $\{r\geq 4\}$. The midpoint
selection
\[
	g_{\mathrm{mid}}(x)=\frac{\mu(x)+1-\nu(x)}{2}=
	\begin{cases}
		1, & r\leq 1, \\[1mm]
		\frac{5}{6}, & 1<r<2, \\[1mm]
		\frac{2}{3}, & r=2, \\[1mm]
		\frac{1}{2}, & 2<r<3, \\[1mm]
		\frac{1}{3}, & r=3, \\[1mm]
		\frac{1}{6}, & 3<r<4, \\[1mm]
		0, & r\geq 4,
	\end{cases}
\]
is discontinuous on each of the four circles $r=1,2,3,4$; on the circles $r=2$
and $r=3$ it takes a value that differs from both one-sided limits. This is
precisely the situation in which $\mu$ and $\nu$ jump on the same set without
compensating each other, so that the existence of a continuous selection is
not evident a priori; it is guaranteed by \Cref{prop:existence}.

The functions \eqref{eq:voting-urysohn} are easily computed. For instance, for
the first expert $\operatorname{dist}(x,A_{1})=(r-2)^{+}$ and
$\operatorname{dist}(x,R_{1})=(3-r)^{+}$, where $t^{+}=\max\{t,0\}$, so that
$u_{1}(x)=\min\{1,(3-r)^{+}\}$. In the same way, $u_{2}(x)=\min\{1,(2-r)^{+}\}$
and $u_{3}(x)=\min\{1,(4-r)^{+}\}$. Each $u_{k}$ decreases linearly from $1$
to $0$ across the annulus on which the $k$-th expert abstains. Averaging,
for $1\leq r\leq 2$ we obtain $\frac{1}{3}\bigl(1+(2-r)+1\bigr)=\frac{4-r}{3}$,
and the same expression is obtained for $2\leq r\leq 3$ and $3\leq r\leq 4$.
Hence
\begin{equation}
	\label{eq:voting-gamma}
	g_{P}(x)=\gamma(x)=\min\Bigl\{1,\ \max\Bigl\{0,\ \frac{4-|x|}{3}\Bigr\}\Bigr\},\qquad
	x\in\mathbb{R}^{2}.
\end{equation}
One checks directly that $\mu\leq\gamma\leq 1-\nu$: for instance, for
$1<r\leq 2$ we have $\gamma(x)\in[\frac{2}{3},1)$, $\mu(x)=\frac{2}{3}$, and
$1-\nu(x)$ equals $1$ for $r<2$ and $\frac{2}{3}$ for $r=2$; the other annuli
are treated in the same way. The consensus selection interpolates linearly
the values imposed by the panel wherever the panel is complete, namely
$1$ for $r\leq 1$, $\frac{2}{3}$ for $r=2$, $\frac{1}{3}$ for $r=3$ and $0$ for
$r\geq 4$.

The consensus selection is also the unique selection of $A$ with the smallest
possible Lipschitz constant, see \Cref{ex:voting-lipschitz}. In terms of
Atanassov's operator $D_{\vartheta}$ (see the \Cref{sec:introduction}), $\gamma$ is the
membership function of $D_{\vartheta}(A)$ with
\[
	\vartheta(x)=2-r\ \ (1<r<2),\qquad \vartheta(x)=3-r\ \ (2<r<3),\qquad
	\vartheta(x)=4-r\ \ (3<r<4),
\]
which jumps from $0$ to $1$ on the circles $r=2$ and $r=3$: the continuous
selection $\gamma$ corresponds to a discontinuous parameter.

\begin{figure}[ht]
	\centering
	\begin{tikzpicture}
		\begin{axis}[
			width=12cm,
			height=7cm,
			xmin=0, xmax=5,
			ymin=-0.05, ymax=1.1,
			axis lines=left,
			xlabel={$r=|x|$ (km)},
			xtick={0,1,2,3,4,5},
			ytick={0,0.3333,0.6667,1},
			yticklabels={$0$,$\frac{1}{3}$,$\frac{2}{3}$,$1$},
			grid=major,
			legend cell align=left,
			legend columns=6,
			legend style={font=\small, at={(0.5,-0.22)}, anchor=north, draw=none},
		]
			% admissible band between mu and 1-nu
			\addplot[name path=lowband, draw=none, const plot, forget plot]
				coordinates {(0,1) (1,0.6667) (2,0.3333) (3,0) (5,0)};
			\addplot[name path=upband, draw=none, const plot, forget plot]
				coordinates {(0,1) (2,0.6667) (3,0.3333) (4,0) (5,0)};
			\addplot[gray!25, forget plot] fill between[of=lowband and upband];
			% mu
			\addplot[blue, very thick, jump mark left, forget plot]
				coordinates {(0,1) (1,0.6667) (2,0.3333) (3,0) (5,0)};
			% 1 - nu
			\addplot[red, very thick, densely dashed, jump mark left, forget plot]
				coordinates {(0,1) (2,0.6667) (3,0.3333) (4,0) (5,0)};
			% midpoint, with its isolated values at r = 2 and r = 3
			\addplot[orange!85!black, very thick, dotted, jump mark left, forget plot]
				coordinates {(0,1) (1,0.8333) (2,0.5) (3,0.1667) (4,0) (5,0)};
			\addplot[only marks, mark=*, mark size=1.8pt, orange!85!black, forget plot]
				coordinates {(2,0.6667) (3,0.3333)};
			% second selection sigma
			\addplot[teal, thick, dashdotted, forget plot, samples=60, domain=1:2]
				{2/3 + (2-x)^2/3};
			\addplot[teal, thick, dashdotted, forget plot, samples=60, domain=2:3]
				{1/3 + (3-x)^2/3};
			\addplot[teal, thick, dashdotted, forget plot, samples=60, domain=3:4]
				{(4-x)^2/3};
			% consensus selection
			\addplot[violet, very thick, forget plot]
				coordinates {(0,1) (1,1) (4,0) (5,0)};
			% legend
			\addlegendimage{area legend, fill=gray!25, draw=none}
			\addlegendentry{$F(x)$}
			\addlegendimage{blue, very thick}
			\addlegendentry{$\mu$}
			\addlegendimage{red, very thick, densely dashed}
			\addlegendentry{$1-\nu$}
			\addlegendimage{orange!85!black, very thick, dotted}
			\addlegendentry{$g_{\mathrm{mid}}$}
			\addlegendimage{violet, very thick}
			\addlegendentry{$\gamma=g_{P}$}
			\addlegendimage{teal, thick, dashdotted}
			\addlegendentry{$\sigma$}
		\end{axis}
	\end{tikzpicture}
	\caption{Radial profiles of the expert-panel intuitionistic fuzzy set of
	\Cref{subsec:voting-example}. The shaded band is the admissible region
	$F(x)$; it collapses to a point for $r\leq 1$, $r=2$, $r=3$ and $r\geq 4$,
	where no expert abstains. The midpoint selection $g_{\mathrm{mid}}$ is
	discontinuous on four circles (the isolated values at $r=2$ and $r=3$ are
	marked). The consensus selection $\gamma=g_{P}$ and the selection $\sigma$
	are two distinct continuous selections.}
	\label{fig:voting}
\end{figure}
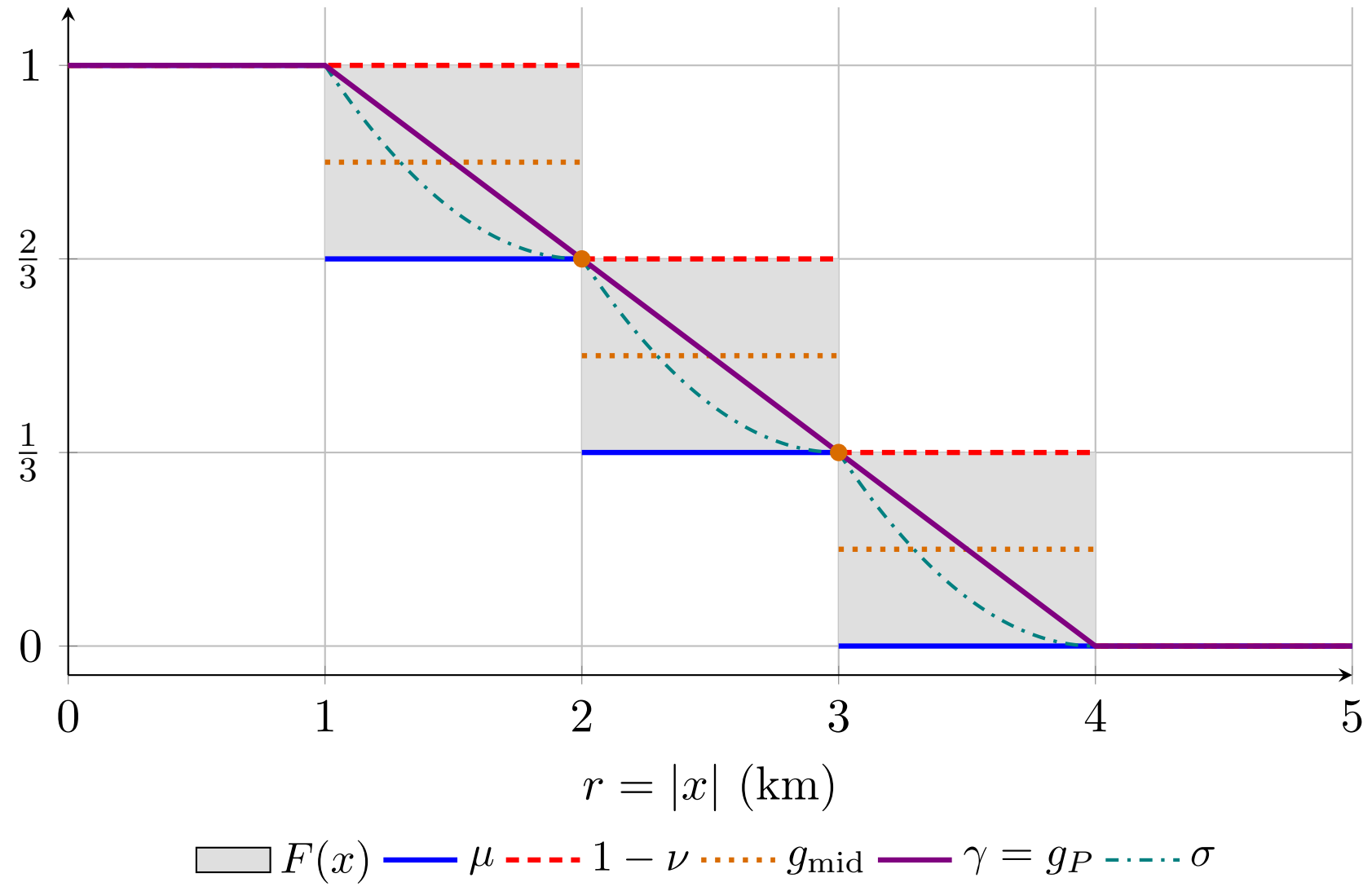

In accordance with the discussion following \Cref{Example1}, the selection is not
unique. A function of the form $g(x)=h(|x|)$ is a continuous selection if and
only if $h\colon[0,\infty)\to\mathbb{R}$ is continuous and satisfies
\[
	h=1\ \text{on }[0,1],\quad h=0\ \text{on }[4,\infty),\quad
	h(2)=\tfrac{2}{3},\quad h(3)=\tfrac{1}{3},
\]
\[
	h(r)\in\bigl[\tfrac{2}{3},1\bigr]\ (1<r<2),\quad
	h(r)\in\bigl[\tfrac{1}{3},\tfrac{2}{3}\bigr]\ (2<r<3),\quad
	h(r)\in\bigl[0,\tfrac{1}{3}\bigr]\ (3<r<4).
\]
Besides $\gamma$, an example is $\sigma(x)=s(|x|)$, where
\[
	s(r)=
	\begin{cases}
		1, & r\leq 1, \\[1mm]
		\frac{2}{3}+\frac{1}{3}(2-r)^{2}, & 1\leq r\leq 2, \\[1mm]
		\frac{1}{3}+\frac{1}{3}(3-r)^{2}, & 2\leq r\leq 3, \\[1mm]
		\frac{1}{3}(4-r)^{2}, & 3\leq r\leq 4, \\[1mm]
		0, & r\geq 4.
	\end{cases}
\]
Selections which are not radial exist as well. All the continuous selections
agree on the set $\{\pi=0\}$, that is, on the disc $\{r\leq 1\}$, on the two
circles $r=2$ and $r=3$, and outside the disc $\{r<4\}$, and they may differ
only on the three open annuli where one expert abstains, the admissible
oscillation being exactly $\pi(x)=\frac{1}{3}$. The consensus selection
$\gamma$ is singled out among them by the structure of the panel, and not by
the topological data alone. The radial profiles of $\mu$, $1-\nu$,
$g_{\mathrm{mid}}$, $\gamma$ and $\sigma$ are represented in \Cref{fig:voting}.

	The radial symmetry of the example has been chosen only to make the
	computations explicit. \Cref{prop:voting-IFS,prop:voting-urysohn} hold for
	arbitrary closed acceptance and rejection regions in an arbitrary metric
	space, for instance non-convex regions bounded by administrative borders, or
	regions in a space of alternatives endowed with a dissimilarity metric. No
	compactness of $X$ has been used. On the other hand, if $X$ is finite and
	endowed with the discrete topology, as in many multi-criteria decision-making
	problems, every function on $X$ is continuous and the selection problem is
	trivial: the results of this paper concern universes that carry a genuine
	topological or metric structure, such as time, space, or measured physical
	quantities.

    The voting interpretation of Atanassov \cite{Ata2} concerns a single decision on each
    element of the universe, the universe being in general finite and not endowed with
    any topology. The model of this section extends that interpretation to a metric
    universe, in which each expert divides the universe into closed regions of
    acceptance and rejection. This spatial structure is what makes upper
    semi-continuity and the explicit selection of \Cref{prop:voting-urysohn}
    available.

\section{Lipschitz selections and stability}
\label{sec:lipschitz}

\subsection{Lipschitz selections}
\label{subsec:lipschitz}

    The existence results of \Cref{sec:cont_def} are qualitative and give no
    information on how rapidly a selection may vary. When $(X,d)$ is a metric
    space, the following result characterises the intuitionistic fuzzy sets that
    admit a selection with a prescribed Lipschitz constant and provides the
    extremal ones explicitly. The construction is an adaptation of McShane--Whitney
    extension formula, see McShane \cite{McS} and Whitney \cite{Whi}. 
    No semi-continuity of $\mu$ and $\nu$
    is required.

    \begin{theorem}
	\label{thm:lipschitz}
	Let $(X,d)$ be a metric space, let $A=(\mu,\nu)$ be an intuitionistic
	fuzzy set in $X$, and let $L\geq 0$. Then there exists an $L$-Lipschitz
	function $g\colon X\to\mathbb{R}$ such that $\mu\leq g\leq 1-\nu$ on $X$ if
	and only if
	\begin{equation}
		\label{eq:lipschitz-condition}
		\mu(x)+\nu(y)\leq 1+L\,d(x,y)\qquad\text{for all }x,y\in X.
	\end{equation}
	In this case, the functions
	\[
		g^{+}_{L}(x)=\inf_{y\in X}\bigl[\,1-\nu(y)+L\,d(x,y)\,\bigr],\qquad
		g^{-}_{L}(x)=\sup_{y\in X}\bigl[\,\mu(y)-L\,d(x,y)\,\bigr]
	\]
	are $L$-Lipschitz selections of $F(x)=[\mu(x),1-\nu(x)]$, and every
	$L$-Lipschitz selection $g$ of $F$ satisfies 
    \[
    g^{-}_{L}\leq g\leq g^{+}_{L}.
    \]
\end{theorem}

\begin{proof}
	\emph{Necessity.} Let $g$ be an $L$-Lipschitz selection. For all $x,y\in X$,
	\[
		\mu(x)-\bigl(1-\nu(y)\bigr)\leq g(x)-g(y)\leq L\,d(x,y),
	\]
	which is \eqref{eq:lipschitz-condition}.

	\emph{Sufficiency.} Assume \eqref{eq:lipschitz-condition}. For every fixed
	$y\in X$, consider the function 
    \[
    \varphi_{y}(x)=1-\nu(y)+L\,d(x,y)
    \]
    of the
	variable $x$. The term $1-\nu(y)$ does not depend on $x$, so no regularity
	of $\nu$ is involved, and by the triangle inequality
	\[
		|\varphi_{y}(x)-\varphi_{y}(x')|=L\,|d(x,y)-d(x',y)|\leq L\,d(x,x')
		\qquad\text{for all }x,x'\in X,
	\]
	that is, each $\varphi_{y}$ is $L$-Lipschitz. Moreover,
	\eqref{eq:lipschitz-condition} gives $\varphi_{y}(x)\geq\mu(x)\geq 0$ for all
	$x,y\in X$, so that $g^{+}_{L}=\inf\limits_{y}\varphi_{y}$ is finite and
	$g^{+}_{L}\geq\mu$. The pointwise infimum of a family of $L$-Lipschitz
	functions which is finite everywhere is $L$-Lipschitz: for all
	$x,x',y\in X$ we have
    \[
    g^{+}_{L}(x)\leq\varphi_{y}(x)\leq\varphi_{y}(x')+L\,d(x,x'),
    \]
	  and taking the
	infimum over $y$ gives $g^{+}_{L}(x)\leq g^{+}_{L}(x')+L\,d(x,x')$. Exchanging
	$x$ and $x'$ yields $|g^{+}_{L}(x)-g^{+}_{L}(x')|\leq L\,d(x,x')$. 
    Choosing $y=x$ in the definition gives
	$g^{+}_{L}(x)\leq 1-\nu(x)$. Thus $g^{+}_{L}$ is an $L$-Lipschitz selection.
	The argument for $g^{-}_{L}$ is symmetric: each function
	$x\mapsto\mu(y)-L\,d(x,y)$ is $L$-Lipschitz and bounded from above by
	$1-\nu(x)$, by \eqref{eq:lipschitz-condition} with the roles of $x$ and $y$
	exchanged, and the choice $y=x$ gives $g^{-}_{L}\geq\mu$.

	\emph{Extremality.} If $g$ is an $L$-Lipschitz selection, then for all
	$x,y\in X$
	\[
		\mu(y)-L\,d(x,y)\leq g(y)-L\,d(x,y)\leq g(x)\leq g(y)+L\,d(x,y)\leq 1-\nu(y)+L\,d(x,y),
	\]
	and taking the supremum and the infimum over $y$ yields
	$g^{-}_{L}(x)\leq g(x)\leq g^{+}_{L}(x)$.
\end{proof}

    \Cref{thm:lipschitz} may be regarded as a two-sided version of the
    McShane--Whitney construction: the extremality of the selections $g^{\pm}_{L}$ is
    the analogue of the corresponding property of the McShane--Whitney extensions,
    proposed by Czipszer--Geh\'er \cite{CziGeh} (see also
    Cobza\c{s}--Miculescu--Nicolae \cite[Theorem 4.1.1]{CobMicNic} and Gutev
    \cite[Theorem 1.2]{Gut}). In particular, \Cref{thm:lipschitz} contains the
    McShane--Whitney extension theorem for $[0,1]$-valued functions, see
    \cite[Theorem 1.3]{Gut}: given $A_{0}\subseteq X$ and an $L$-Lipschitz function
    $\varphi\colon A_{0}\to[0,1]$, let $\mu=\varphi$ and $\nu=1-\varphi$ on $A_{0}$, and
    $\mu=\nu=0$ on $X\setminus A_{0}$. The selections of $F$ are then exactly the
    extensions of $\varphi$ with values in $[0,1]$, and condition
    \eqref{eq:lipschitz-condition} below reduces to the $L$-Lipschitz condition for
    $\varphi$ on $A_{0}$.
    
    For $x=y$, condition \eqref{eq:lipschitz-condition} reduces to the defining
    inequality $\mu(x)+\nu(x)\leq 1$ of an intuitionistic fuzzy set. It may therefore
    be read as a non-local version of that inequality, in which the membership
    degree at one point and the non-membership degree at another point may exceed
    the total amount $1$ by at most $L$ times their distance. The smallest constant
    for which \eqref{eq:lipschitz-condition} holds,
    \[
    	L^{*}(A)=\sup_{x\neq y}\frac{\bigl(\mu(x)+\nu(y)-1\bigr)^{+}}{d(x,y)}\in[0,+\infty],
    \]
    is thus the smallest Lipschitz constant of a selection of $F$, whenever it is
    finite: indeed, by the definition of supremum, \eqref{eq:lipschitz-condition}
    holds with $L=L^{*}(A)$. In particular, $L^{*}(A)=0$ if and only if
    $\sup\limits_{X}\mu+\sup\limits_{X}\nu\leq 1$, that is, if and only if $F$ admits a constant
    selection.

\begin{remark}
	\label{rem:independence}
	On a metric space, \Cref{thm:lipschitz} provides continuous selections without
	any semi-continuity assumption, whenever $L^{*}(A)<\infty$. The two sufficient
	conditions for the existence of a continuous selection, namely the upper
	semi-continuity of $\mu$ and $\nu$ (\Cref{prop:existence}) and the finiteness of
	$L^{*}(A)$ (\Cref{thm:lipschitz}), are independent. In
	\Cref{rem:voting-non-uniform}, $\mu$ and $\nu$ are u.s.c.\ while
	$L^{*}(A)=+\infty$. Conversely, let $X=[0,1]$, $\mu=\frac{1}{2}\mathbf{1}_{\mathbb{Q}\cap[0,1]}$
	and $\nu\equiv 0$: the function $\mu$ is not u.s.c., since
	$\{x:\mu(x)\geq\frac{1}{2}\}=\mathbb{Q}\cap[0,1]$ is not closed, while
	$\mu(x)+\nu(y)\leq\frac{1}{2}$ for all $x,y$, so that $L^{*}(A)=0$ and the constant
	$\frac{1}{2}$ is a selection.
\end{remark}

\begin{example}
	\label{ex:voting-lipschitz}
	Consider the intuitionistic fuzzy set of \Cref{subsec:voting-example}. We claim
	that $L^{*}(A)=\frac{1}{3}$ and that the consensus selection $\gamma$ in
	\eqref{eq:voting-gamma} is the unique selection with Lipschitz constant
	$\frac{1}{3}$. The function $\gamma$ is the composition of the $1$-Lipschitz map
	$x\mapsto|x|$ with a piecewise linear function of slope at most $\frac{1}{3}$ in
	absolute value, hence it is $\frac{1}{3}$-Lipschitz, and it is a selection; by
	\Cref{thm:lipschitz}, $L^{*}(A)\leq\frac{1}{3}$. Conversely, let $g$ be an
	$L$-Lipschitz selection and fix a unit vector $e\in\mathbb{R}^{2}$. Since $\pi=0$
	on the circles $r=1$ and $r=4$, the discussion at the end of \Cref{sec:cont_def} gives $g(e)=1$ and
	$g(4e)=0$, whence $1\leq 3L$. If $L=\frac{1}{3}$ and $1\leq r\leq 4$, then
	\[
		g(re)\leq g(4e)+\tfrac{1}{3}(4-r)=\tfrac{4-r}{3},\qquad
		g(re)\geq g(e)-\tfrac{1}{3}(r-1)=\tfrac{4-r}{3},
	\]
	so that $g(re)=\gamma(re)$. Since $g=1$ on $\{r\leq 1\}$ and $g=0$ on
	$\{r\geq 4\}$, we conclude that $g=\gamma$. In particular
	$g^{-}_{1/3}=g^{+}_{1/3}=\gamma$. For larger constants uniqueness is lost: the
	selection $\sigma$ of \Cref{subsec:voting-example} is $\frac{2}{3}$-Lipschitz, so
	that for $L\geq\frac{2}{3}$ both $\gamma$ and $\sigma$ are $L$-Lipschitz
	selections.
\end{example}

\subsection{Locally Lipschitz selections}
\label{subsec:loc-lipschitz}

Between the continuous selections of \Cref{prop:existence} and the Lipschitz
selections of \Cref{thm:lipschitz} there is an intermediate class. Recall that a
function $g\colon X\to\mathbb{R}$ on a metric space $(X,d)$ is \emph{locally
Lipschitz} if every point of $X$ has a neighbourhood on which $g$ is Lipschitz.
For an intuitionistic fuzzy set $A=(\mu,\nu)$ in $X$ and a set $U\subseteq X$ we
write
\[
	L^{*}(A;U)=\sup_{x,y\in U,\ x\neq y}\frac{\bigl(\mu(x)+\nu(y)-1\bigr)^{+}}{d(x,y)}\in[0,+\infty],
\]
so that $L^{*}(A;X)=L^{*}(A)$. The following result is the local counterpart of
\Cref{thm:lipschitz}; it relies on the characterisation of locally Lipschitz
functions and on the Lipschitz partitions of unity discussed by Gutev \cite{Gut}.

\begin{theorem}
	\label{thm:loc-lipschitz}
	Let $(X,d)$ be a metric space and let $A=(\mu,\nu)$ be an intuitionistic fuzzy
	set in $X$. The following conditions are equivalent:
	\begin{itemize}
		\item[\textnormal{(a)}] $F(x)=[\mu(x),1-\nu(x)]$ admits a locally Lipschitz
			selection;
		\item[\textnormal{(b)}] every point of $X$ has an open neighbourhood $U$ such
			that $L^{*}(A;U)<\infty$;
		\item[\textnormal{(c)}] there exists a continuous function
			$\Lambda\colon X\times X\to[0,+\infty)$ such that
			\begin{equation}
				\label{eq:loc-lipschitz-condition}
				\mu(x)+\nu(y)\leq 1+\Lambda(x,y)\,d(x,y)\qquad\text{for all }x,y\in X.
			\end{equation}
	\end{itemize}
\end{theorem}

\begin{proof}
	(a)$\Rightarrow$(c). Let $g$ be a locally Lipschitz selection of $F$. By
	\cite[Theorem 6.5]{Gut} there exists a continuous function
	$\Lambda\colon X\times X\to[0,+\infty)$ such that
	$|g(x)-g(y)|\leq\Lambda(x,y)\,d(x,y)$ for all $x,y\in X$. Hence
	\[
		\mu(x)-\bigl(1-\nu(y)\bigr)\leq g(x)-g(y)\leq\Lambda(x,y)\,d(x,y),
	\]
	which is \eqref{eq:loc-lipschitz-condition}.

	(c)$\Rightarrow$(b). Let $p\in X$. Since $\Lambda$ is continuous at $(p,p)$, there
	exists $r>0$ such that $\Lambda\leq\Lambda(p,p)+1$ on $U\times U$, where
	$U=\{x\in X:\ d(x,p)<r\}$. By \eqref{eq:loc-lipschitz-condition},
	$L^{*}(A;U)\leq\Lambda(p,p)+1$.

	(b)$\Rightarrow$(a). For every $p\in X$ choose an open neighbourhood $U_{p}$ with
	$L_{p}=L^{*}(A;U_{p})<\infty$. Applying \Cref{thm:lipschitz} to the restriction of
	$A$ to the metric space $(U_{p},d)$, we obtain an $L_{p}$-Lipschitz function
	$g_{p}\colon U_{p}\to[0,1]$ with $\mu\leq g_{p}\leq 1-\nu$ on $U_{p}$. By a theorem
	of Frol\'{\i}k \cite[Theorem 1]{Fro} (see also \cite[Section 5]{Gut}), 
    we know that each open cover of a metric space admits a locally finite partition of unity consisting of Lipschitz functions, then 
    the open
	cover $\{U_{p}:\ p\in X\}$ admits a locally finite partition of unity
	$\{\xi_{\alpha}:\ \alpha\in\mathcal{A}\}$ consisting of Lipschitz functions and
	subordinated to it: for every $\alpha\in\mathcal{A}$ there is $p(\alpha)\in X$
	with $\operatorname{supp}(\xi_{\alpha})\subseteq U_{p(\alpha)}$. Define
	$h_{\alpha}=\xi_{\alpha}\,g_{p(\alpha)}$ on $U_{p(\alpha)}$ and $h_{\alpha}=0$ on
	$X\setminus U_{p(\alpha)}$. On the open set $U_{p(\alpha)}$, $h_{\alpha}$ is the
	product of two bounded Lipschitz functions, hence it is Lipschitz. On the open
	set $X\setminus\operatorname{supp}(\xi_{\alpha})$ it vanishes. Since these two open
	sets cover $X$, $h_{\alpha}$ is locally Lipschitz. The function
	$g=\sum_{\alpha}h_{\alpha}$ is a locally finite sum of locally Lipschitz functions,
	and is therefore locally Lipschitz. Finally, for $x\in X$ we have
	\[
		g(x)=\sum_{\alpha:\ \xi_{\alpha}(x)>0}\xi_{\alpha}(x)\,g_{p(\alpha)}(x),
	\]
	a finite convex combination of points of $F(x)$, because $\xi_{\alpha}(x)>0$
	implies $x\in U_{p(\alpha)}$ and hence $g_{p(\alpha)}(x)\in F(x)$. Since $F(x)$ is an
	interval, $g(x)\in F(x)$.
\end{proof}

Condition \eqref{eq:loc-lipschitz-condition} is a non-uniform version of
\eqref{eq:lipschitz-condition}, in which the constant $L$ is replaced by a
continuous function of the pair of points. As before, for $x=y$ it reduces to the
inequality $\mu(x)+\nu(x)\leq 1$. Two consequences of \Cref{thm:loc-lipschitz} are
the following.

\begin{corollary}
	\label{cor:loc-lip-usc}
	Let $(X,d)$ be a metric space and let $A=(\mu,\nu)$ be an intuitionistic fuzzy
	set in $X$ such that $\mu$ and $\nu$ are upper semi-continuous and $\pi(x)>0$ for
	every $x\in X$. Then $F$ admits a locally Lipschitz selection.
\end{corollary}

\begin{proof}
	Let $p\in X$ and $\varepsilon=\pi(p)/3>0$. By upper semi-continuity there exists
	an open neighbourhood $U$ of $p$ such that $\mu(x)<\mu(p)+\varepsilon$ and
	$\nu(y)<\nu(p)+\varepsilon$ for all $x,y\in U$. Then
	$\mu(x)+\nu(y)<1-\pi(p)+2\varepsilon<1$ for all $x,y\in U$, so that
	$L^{*}(A;U)=0$, and the conclusion follows from \Cref{thm:loc-lipschitz}.
\end{proof}

\Cref{cor:loc-lip-usc} can also be deduced from the locally Lipschitz version of
Dowker's insertion theorem obtained by Gutev \cite[Corollary 7.5]{Gut}, which
moreover provides a selection $g$ with $\mu<g<1-\nu$.

\begin{corollary}
	\label{cor:voting-loc-lip}
	In the voting model of \Cref{sec:voting}, the consensus selection $g_{P}$ is
	locally Lipschitz. If, moreover,
	$\delta_{k}=\operatorname{dist}(A_{k},R_{k})>0$ for every $k$ (with $\delta_{k}=+\infty$
	if $A_{k}$ or $R_{k}$ is empty), then $g_{P}$ is Lipschitz with constant at most
	$\frac{3}{n}\sum_{k=1}^{n}\delta_{k}^{-1}$.
\end{corollary}

\begin{proof}
	It suffices to consider an expert $k$ with $A_{k}\neq\emptyset\neq R_{k}$, since
	otherwise $u_{k}$ is constant. Write $u_{k}=f/s$ with
	$f=\operatorname{dist}(\cdot,R_{k})$ and $s=\operatorname{dist}(\cdot,A_{k})+f$. The function
	$f$ is $1$-Lipschitz, $s$ is $2$-Lipschitz, $0\leq f\leq s$ and $s>0$ on $X$. For
	$x,y\in X$,
	\[
		|u_{k}(x)-u_{k}(y)|
		=\frac{\bigl|(f(x)-f(y))\,s(y)+f(y)\,(s(y)-s(x))\bigr|}{s(x)\,s(y)}
		\leq\frac{3\,d(x,y)}{s(x)}.
	\]
	Since $s$ is continuous and positive, for every $p\in X$ we have $s\geq s(p)/2$ in a
	neighbourhood of $p$, on which $u_{k}$ is therefore $6/s(p)$-Lipschitz; hence
	$u_{k}$ and $g_{P}$ are locally Lipschitz. Moreover, by the triangle inequality
	$s(x)\geq\delta_{k}$ for every $x\in X$, so that $u_{k}$ is $3/\delta_{k}$-Lipschitz,
	and the bound for $g_{P}$ follows by averaging.
\end{proof}

The three classes of selections considered so far are distinct, as the following
examples show.

\begin{example}
	\label{ex:hierarchy}
    The following are possible selections examples:
	\begin{itemize}
		\item[\textnormal{(i)}] Let $X=[0,1]$, $\mu(x)=\sqrt{x}$ and $\nu(x)=1-\sqrt{x}$.
			Then $\mu$ and $\nu$ are continuous and $\pi\equiv 0$, so that 
            the discussion at the end of \Cref{sec:cont_def}
			 ensures that the unique selection is $g(x)=\sqrt{x}$, which is
			continuous but not locally Lipschitz at $0$. Accordingly, condition (b) of
			\Cref{thm:loc-lipschitz} fails at $0$: for every neighbourhood $U$ of $0$ and
			$x\in U\setminus\{0\}$, we have
			$\bigl(\mu(x)+\nu(0)-1\bigr)/x=1/\sqrt{x}$, which is unbounded.
		\item[\textnormal{(ii)}] In the example of \Cref{rem:voting-non-uniform}, the
			consensus selection, that is, the Urysohn function of $(A_{1},R_{1})$, is
			locally Lipschitz by \Cref{cor:voting-loc-lip}, while no selection is
			Lipschitz, since $L^{*}(A)=+\infty$.
		\item[\textnormal{(iii)}] In \Cref{ex:voting-lipschitz} the consensus selection
			is Lipschitz, with the smallest possible constant $L^{*}(A)=\frac{1}{3}$.
	\end{itemize}
	Thus the intuitionistic fuzzy sets admitting a Lipschitz selection form a proper
	subclass of those admitting a locally Lipschitz selection, which in turn form a
	proper subclass of those admitting a continuous selection.
\end{example}

\subsection{Stability with respect to the Hausdorff distance}
\label{subsec:stability}

In the applications discussed in \Cref{sec:introduction}, the functions $\mu$
and $\nu$ are obtained from expert judgements or measurements and are therefore
known only up to some error. It is then natural to ask how the family of
continuous selections changes when the intuitionistic fuzzy set is perturbed.
To measure the perturbation we use the distance introduced by Grzegorzewski
\cite{Grz}, which is based on the Hausdorff distance between the intervals
of the form
$F_{A}(x)=[\mu_{A}(x),1-\nu_{A}(x)]$. For two closed intervals of $\mathbb{R}$,
the Hausdorff distance is the largest distance between corresponding
endpoints \cite[formula (18)]{Grz}. Hence, for intuitionistic fuzzy sets
$A=(\mu_{A},\nu_{A})$ and $B=(\mu_{B},\nu_{B})$ in $X$,
\[
	d_{H}\bigl(F_{A}(x),F_{B}(x)\bigr)
	=\max\bigl\{|\mu_{A}(x)-\mu_{B}(x)|,\ |\nu_{A}(x)-\nu_{B}(x)|\bigr\},
	\qquad x\in X.
\]
Grzegorzewski aggregates these pointwise distances by a finite sum, since the
universe considered in \cite{Grz} is finite. On an arbitrary universe we take
the supremum and set
\begin{equation}
	\label{eq:dist-infty}
	d_{\infty}(A,B)=\sup_{x\in X}\max\bigl\{|\mu_{A}(x)-\mu_{B}(x)|,\
	|\nu_{A}(x)-\nu_{B}(x)|\bigr\}\in[0,1].
\end{equation}
It is immediate that $d_{\infty}$ is a metric on the family of
intuitionistic fuzzy sets in $X$, and that it reduces to the uniform distance
$\sup_{X}|\mu_{A}-\mu_{B}|$ when $A$ and $B$ are fuzzy sets, that is, when
$\nu_{A}=1-\mu_{A}$ and $\nu_{B}=1-\mu_{B}$.

Let $C_{b}(X)$ denote the space of bounded continuous real functions on $X$
with the norm $\Vert g\Vert_{\infty}=\sup\limits_{X}|g|$, and denote by
\[
	\mathcal{S}(A)=\bigl\{g\in C_{b}(X):\ \mu_{A}\leq g\leq 1-\nu_{A}\ \text{on }X\bigr\}
\]
the set of continuous selections of $F_{A}$. Every continuous selection takes
values in $[0,1]$ and therefore belongs to $C_{b}(X)$. For nonempty subsets
$\mathcal{U},\mathcal{V}\subseteq C_{b}(X)$ we write
\[
	D_{H}(\mathcal{U},\mathcal{V})=\max\Bigl\{\sup_{g\in\mathcal{U}}\inf_{h\in\mathcal{V}}\Vert g-h\Vert_{\infty},\
	\sup_{h\in\mathcal{V}}\inf_{g\in\mathcal{U}}\Vert g-h\Vert_{\infty}\Bigr\}
\]
for their Hausdorff distance in $C_{b}(X)$.

\begin{theorem}
	\label{thm:stability}
	Let $X$ be a $T_{4}$ space, and let $A=(\mu_{A},\nu_{A})$ and
	$B=(\mu_{B},\nu_{B})$ be intuitionistic fuzzy sets in $X$ whose membership
	and non-membership functions are upper semi-continuous. Then for every
	$g\in\mathcal{S}(A)$ there exists $h\in\mathcal{S}(B)$ such that
	\[
		\Vert g-h\Vert_{\infty}\leq d_{\infty}(A,B).
	\]
	Consequently,
	\begin{equation}
		\label{eq:stability}
		D_{H}\bigl(\mathcal{S}(A),\mathcal{S}(B)\bigr)\leq d_{\infty}(A,B).
	\end{equation}
\end{theorem}

\begin{proof}
	By \Cref{prop:existence}, the sets $\mathcal{S}(A)$ and $\mathcal{S}(B)$ are
	nonempty. Let $\delta=d_{\infty}(A,B)$ and $g\in\mathcal{S}(A)$, and define
	\[
		a(x)=\max\bigl\{\mu_{B}(x),\ g(x)-\delta\bigr\},\qquad
		b(x)=\min\bigl\{1-\nu_{B}(x),\ g(x)+\delta\bigr\}.
	\]
	The function $a$ is the maximum of two upper semi-continuous functions, and
	is therefore upper semi-continuous; similarly, $b$ is lower semi-continuous,
	since $1-\nu_{B}$ is lower semi-continuous and $g+\delta$ is continuous. We
	claim that $a\leq b$ on $X$, which amounts to the four inequalities
	\[
		\mu_{B}\leq 1-\nu_{B},\qquad \mu_{B}\leq g+\delta,\qquad
		g-\delta\leq 1-\nu_{B},\qquad g-\delta\leq g+\delta.
	\]
	The first one holds because $B$ is an intuitionistic fuzzy set and the last
	one is trivial. By \eqref{eq:dist-infty} and since $g\in\mathcal{S}(A)$, we have
	\[
		\mu_{B}\leq\mu_{A}+\delta\leq g+\delta,\qquad
		g-\delta\leq 1-\nu_{A}-\delta\leq 1-\nu_{B},
	\]
	which are the second and the third inequality. By the Kat\v{e}tov--Tong
	theorem (\Cref{thm:KatTon}) there exists a continuous function
	$h\colon X\to\mathbb{R}$ with $a\leq h\leq b$. Then
	$\mu_{B}\leq h\leq 1-\nu_{B}$, so that $h\in\mathcal{S}(B)$, and
	$g-\delta\leq h\leq g+\delta$, that is, $\Vert g-h\Vert_{\infty}\leq\delta$.

	This proves that
	$\sup\limits_{g\in\mathcal{S}(A)}\inf_{h\in\mathcal{S}(B)}\Vert g-h\Vert_{\infty}\leq\delta$.
	Exchanging the roles of $A$ and $B$ gives the same bound for the other term
	in the definition of $D_{H}$, and \eqref{eq:stability} follows.
\end{proof}

The estimate \eqref{eq:stability} cannot be improved in general. If
$\pi_{A}\equiv 0$ and $\pi_{B}\equiv 0$, then again by the discussion at the end of \Cref{sec:cont_def},
$\mathcal{S}(A)=\{\mu_{A}\}$ and $\mathcal{S}(B)=\{\mu_{B}\}$, and both sides of
\eqref{eq:stability} are equal to $\sup\limits_{X}|\mu_{A}-\mu_{B}|$. This is the
counterpart of the fact, observed in \cite[Lemma 5]{Grz}, that the distances
based on the Hausdorff metric reduce to the classical ones for fuzzy sets.

\Cref{thm:stability} admits the following interpretation. Following
\cite[Example 1]{Grz}, the interval $F_{A}(x)$ collects all the possible outcomes
of a negotiation once the hesitation has been removed, and a continuous selection
is a choice of such outcomes that varies continuously over the universe. The
theorem states that, if the expert judgements are modified by at most $\delta$
in the sense of \eqref{eq:dist-infty}, then every such choice can be adapted
to the new judgements by modifying it by at most $\delta$, uniformly on $X$.
Note that the midpoint selection satisfies the analogous estimate
$\Vert g_{\mathrm{mid},A}-g_{\mathrm{mid},B}\Vert_{\infty}\leq d_{\infty}(A,B)$, but
it is not continuous in general; \Cref{thm:stability} provides the same
quantitative control within the class of continuous selections.

When $X$ is a metric space, a version of \Cref{thm:stability} holds for
Lipschitz selections, without any semi-continuity assumption and with an
explicit choice of $h$. For $L\geq 0$ let $\mathcal{S}_{L}(A)$ denote the set
of $L$-Lipschitz selections of $F_{A}$, and let $g^{\pm}_{L,A}$ be the
extremal selections given by \Cref{thm:lipschitz}.

\begin{proposition}
	\label{prop:stability-lipschitz}
	Let $(X,d)$ be a metric space, let $A$ and $B$ be intuitionistic fuzzy sets
	in $X$, and let $L\geq 0$ be such that $\mathcal{S}_{L}(A)$ and
	$\mathcal{S}_{L}(B)$ are nonempty. Then
	\[
		\Vert g^{+}_{L,A}-g^{+}_{L,B}\Vert_{\infty}\leq d_{\infty}(A,B),\qquad
		\Vert g^{-}_{L,A}-g^{-}_{L,B}\Vert_{\infty}\leq d_{\infty}(A,B),
	\]
	and, for every $g\in\mathcal{S}_{L}(A)$, the function
	\[
		h=\min\bigl\{\max\{g,\ g^{-}_{L,B}\},\ g^{+}_{L,B}\bigr\}
	\]
	belongs to $\mathcal{S}_{L}(B)$ and satisfies
	$\Vert g-h\Vert_{\infty}\leq d_{\infty}(A,B)$. In particular,
	$D_{H}\bigl(\mathcal{S}_{L}(A),\mathcal{S}_{L}(B)\bigr)\leq d_{\infty}(A,B)$.
\end{proposition}

\begin{proof}
	Let $\delta=d_{\infty}(A,B)$. For all $x,y\in X$ we have
    \[
    1-\nu_{B}(y)+L\,d(x,y)\leq 1-\nu_{A}(y)+L\,d(x,y)+\delta.
    \]
	Taking the
	infimum over $y$ gives $g^{+}_{L,B}\leq g^{+}_{L,A}+\delta$, and exchanging
	$A$ and $B$ yields the first estimate. The second one is obtained in the
	same way from 
    \[
    \mu_{B}(y)-L\,d(x,y)\leq\mu_{A}(y)-L\,d(x,y)+\delta.
    \]

	Let $g\in\mathcal{S}_{L}(A)$. Since $\mathcal{S}_{L}(B)$ is nonempty, its
	elements lie between $g^{-}_{L,B}$ and $g^{+}_{L,B}$ by \Cref{thm:lipschitz},
	so that $g^{-}_{L,B}\leq g^{+}_{L,B}$ and $h$ takes values between
	$g^{-}_{L,B}$ and $g^{+}_{L,B}$. Hence $\mu_{B}\leq h\leq 1-\nu_{B}$. Moreover,
	$h$ is $L$-Lipschitz, being obtained from $L$-Lipschitz functions by taking
	pointwise maxima and minima. Thus $h\in\mathcal{S}_{L}(B)$. Finally, at a
	point where $h\neq g$ we have either $g<h=g^{-}_{L,B}$, in which case
	\[
		0<h-g\leq g^{-}_{L,B}-g^{-}_{L,A}\leq\delta,
	\]
	since $g\geq g^{-}_{L,A}$, or $g>h=g^{+}_{L,B}$, in which case
	\[
		0<g-h\leq g^{+}_{L,A}-g^{+}_{L,B}\leq\delta,
	\]
	since $g\leq g^{+}_{L,A}$. Therefore $\Vert g-h\Vert_{\infty}\leq\delta$, and the
	bound on $D_{H}$ follows by exchanging the roles of $A$ and $B$.
\end{proof}

In particular, the maps $A\mapsto g^{\pm}_{L,A}$, which single out an
explicit selection for each intuitionistic fuzzy set, are $1$-Lipschitz from
the metric space $(\mathrm{IFS}(X),d_{\infty})$ into $C_{b}(X)$, on the
family of intuitionistic fuzzy sets that admit $L$-Lipschitz selections.

\section{Relations with other fuzzy sets}
	\label{sec:conclusion} In this section we discuss some applications of continuous
	selections, their relation with interval type-2 fuzzy sets, the extension of the
	present approach to other generalizations of intuitionistic fuzzy sets, and the
	choice of a distinguished selection.

	\subsection{Continuous selections and real-world applications}
	The midpoint selection \eqref{eq:midpoint-introduction}, that is, the membership
	function of $D_{1/2}(A)$, is a widely adopted deterministic choice. When $\mu$ and
	$\nu$ are merely upper semi-continuous, however, their graphs may exhibit jumps, and
	$g_{\mathrm{mid}}$ inherits these jumps whenever $\mu-\nu$ is discontinuous, as in
	\Cref{subsec:voting-example}; in a control loop this produces a discontinuous
	control signal. Some situations in which the continuity of the selected
	membership function matters are the following.
	\begin{itemize}
		\item \textbf{Fuzzy Control under Sensor Uncertainty.} Sensor measurements
			$x$ are often corrupted by noise and can be modelled as intuitionistic
			fuzzy sets (IFSs). The controller is required to compute a control action based
			on these uncertain inputs. The existence of a continuous selection
			function $g(x)$ provides a theoretical guarantee that a robust and smooth
			control law can be designed, preventing discontinuities in the control output
			caused by abrupt variations in sensor measurements.

			In robotics and nonlinear control, discontinuous control signals are a
			known source of the chattering phenomenon, extensively documented in Sliding
			Mode Control. A jump in the control law would require the actuator to
			switch in an arbitrarily short time, which physical actuators cannot do:
			their bandwidth is finite and their moving parts have inertia. The resulting
			high-frequency oscillation excites unmodelled dynamics and increases mechanical
			wear and energy consumption. Continuity of the selected membership
			function is a necessary condition for avoiding this behaviour at the
			defuzzification stage; it is of course not sufficient for the closed-loop system
			as a whole, whose analysis lies outside the scope of this paper. We refer to
			Hasan--Alam \cite{HasAla} for a discussion of chattering in exoskeletons
			for rehabilitation applications.

		\item \textbf{Consensus Aggregation in Group Decision-Making.} Suppose that
			the interval $[\mu(x),\,1-\nu(x)]$ represents the range of agreement and
			disagreement among a group of experts regarding a parameter $x$. In this setting,
			the continuous selection function serves as a mathematical mechanism for
			extracting a compromise score that varies continuously with the experts' opinions,
			thereby avoiding artificial or illogical jumps resulting from marginal
			changes in individual assessments. \Cref{thm:stability} makes this
			quantitative: a perturbation of size $\delta$ of the experts' judgements, in
			the distance \eqref{eq:dist-infty}, moves the set of continuous selections by
			at most $\delta$.
	\end{itemize}

	\subsection{Relation with interval type-2 fuzzy sets}
	The interval-valued map $F(x)=[\mu(x),1-\nu(x)]$ coincides, as a geometric
	object, with the Footprint of Uncertainty (FOU) of an Interval Type-2 Fuzzy Set,
	see Mendel--John \cite{MenJoh}: the membership function $\mu$ plays the role
	of the lower membership function, its counterpart $1-\nu$ that of the upper membership
	function, and the hesitation degree $\pi(x)$ measures the width of the
	footprint at $x$. The correspondence between the two frameworks is well
	documented, see, e.g., Bustince--Burillo \cite{BusBur} and Deschrijver--Kerre
	\cite{DesKer} for the relations between intuitionistic, interval-valued and
	related fuzzy sets.

	The two settings differ in what is extracted from the footprint. Type-reduction
	procedures, such as the Karnik--Mendel algorithm \cite{KarMen}, aggregate the embedded type-1
	sets and return a centroid interval, from which a crisp value is obtained by a
	further step. The construction considered here does not aggregate:
	\Cref{prop:existence} produces one embedded type-1 membership
	function, selected so as to be continuous. The two operations answer different
	questions and their outputs are not directly comparable. What the selection approach
	adds is a continuity guarantee on the resulting membership function, obtained
	from the semi-continuity of $\mu$ and $\nu$ alone, or, on metric spaces, from the
	finiteness of $L^{*}(A)$.

	\subsection{Other generalizations of intuitionistic fuzzy sets}
	\label{subsec:other-sets}
	Several generalizations of intuitionistic fuzzy sets reduce to the present
	setting as far as the selection problem is concerned. Pythagorean fuzzy sets
	\cite{Yag1,Yag2}, defined by the condition $\mu^{2}+\nu^{2}\leq 1$, coincide with
	the intuitionistic fuzzy sets of second type already considered by Atanassov
	\cite{Ata2}, and $q$-rung orthopair fuzzy sets \cite{Yag3}, defined by
	$\mu^{q}+\nu^{q}\leq 1$, are of the same nature. More generally, let
	$\varphi\colon[0,1]\to[0,1]$ be an increasing homeomorphism and let $\mu,\nu$
	satisfy $\varphi(\mu)+\varphi(\nu)\leq 1$. The admissible values of a compatible
	membership function form the interval
	$[\mu(x),\varphi^{-1}(1-\varphi(\nu(x)))]$, and $g$ is a continuous selection of
	this interval-valued map if and only if $\varphi\circ g$ is a continuous selection
	of $[\varphi\circ\mu,1-\varphi\circ\nu]$. Since $\varphi$ is increasing and
	continuous, $\varphi\circ\mu$ and $\varphi\circ\nu$ are u.s.c.\ whenever $\mu$ and
	$\nu$ are, and $(\varphi\circ\mu,\varphi\circ\nu)$ is an intuitionistic fuzzy set.
	Hence \Cref{prop:existence} and the results of \Cref{sec:voting} apply straightforwardly,
	while the quantitative results of \Cref{sec:lipschitz} transfer with constants
	depending on the moduli of continuity of $\varphi$ and $\varphi^{-1}$.

	Similarly, a membership function compatible with a picture fuzzy set
	$(\mu,\eta,\nu)$, see, e.g., \cite{Cuo}, with refusal degree $\varrho=1-\mu-\eta-\nu$, must
	satisfy $\mu\leq g\leq\mu+\eta+\varrho=1-\nu$: the admissible interval coincides
	with that of the intuitionistic fuzzy set $(\mu,\nu)$, and picture fuzzy sets do not
	give rise to a new selection problem. The same holds for spherical fuzzy sets
	\cite{KutKah}, since $\mu^{2}+\eta^{2}+\nu^{2}\leq 1$ implies $\mu^{2}+\nu^{2}\leq 1$,
	and for single-valued neutrosophic sets \cite{Sma,WanSmaZhaSun} whenever the
	truth and falsity degrees satisfy $T+F\leq 1$, the admissible interval being
	$[T,1-F]$.

	Genuinely higher-dimensional selection problems arise when the degrees
	themselves are interval-valued, as for interval-valued intuitionistic fuzzy sets
	\cite{AtaGar}. In this case the Kat\v{e}tov--Tong theorem has to be replaced by
	Michael's theorem (\Cref{Thm:MichaelSelection}), see \Cref{rem:michael-route}, or by
	lattice-valued insertion theorems \cite{KubMar}, and the lower semi-continuity of
	the corresponding set-valued map has to be established.

	\subsection{A criterion for choosing a selection}
	\label{subsec:criterion}
	As shown in \Cref{Example1} and in \Cref{subsec:voting-example}, in general the
	continuous selection is not unique, and a distinguished element can be singled
	out only by an additional criterion. When $\mu$ and $\nu$ are continuous, as for
	triangular, trapezoidal or Gaussian membership functions, every function
	$\mu+\vartheta\pi$ with $\vartheta\colon X\to[0,1]$ continuous is a continuous
	selection, and the midpoint $g_{\mathrm{mid}}$, corresponding to
	$\vartheta\equiv\frac{1}{2}$, is a natural choice. When $\mu$ and $\nu$ are only
	upper semi-continuous, $g_{\mathrm{mid}}$ may be discontinuous, and a different
	criterion is needed.

	The results of \Cref{sec:lipschitz} make such a criterion available. Let $(X,d)$
	be a compact metric space, let $A$ be an intuitionistic fuzzy set in $X$ with
	$\mu$ and $\nu$ Borel measurable (for instance u.s.c.), and let $L\geq L^{*}(A)$ be
	finite. The set $\mathcal{S}_{L}(A)$ of $L$-Lipschitz selections is nonempty by
	\Cref{thm:lipschitz} and convex. Its elements take values in $[0,1]$ and are
	equi-Lipschitz, so that $\mathcal{S}_{L}(A)$ is relatively compact in $C(X)$ by the
	Arzel\`a--Ascoli theorem. It is also closed in $C(X)$, since the inequalities
	$\mu\leq g\leq 1-\nu$ and $|g(x)-g(y)|\leq L\,d(x,y)$ are preserved under uniform
	convergence. Hence $\mathcal{S}_{L}(A)$ is a compact convex subset of $C(X)$. Let
	$m$ be a finite Borel measure on $X$ with full support, for instance the Lebesgue
	measure on a compact interval, and consider
	\[
		J(g)=\int_{X}\bigl(g-g_{\mathrm{mid}}\bigr)^{2}\,dm .
	\]
	The functional $J$ is continuous on $C(X)$ and strictly convex on
	$\mathcal{S}_{L}(A)$: if $g_{1}\neq g_{2}$ are continuous, then $g_{1}-g_{2}$ does not
	vanish on a nonempty open set, which has positive measure, so that
	\[
		\int_{X}(g_{1}-g_{2})^{2}\,dm>0 .
	\]
	Therefore $J$ has a unique minimizer $g^{*}_{L}$ on $\mathcal{S}_{L}(A)$, the
	$L$-Lipschitz selection which is closest to the midpoint in the mean-square
	sense. When $L=L^{*}(A)$ the set $\mathcal{S}_{L}(A)$ may already reduce to a single
	element, as in \Cref{ex:voting-lipschitz}.

	This criterion is also related to the Hausdorff distance used by Grzegorzewski
	\cite{Grz}: for $g(x)\in F(x)$ one has
	\[
		d_{H}\bigl(\{g(x)\},F(x)\bigr)=\frac{\pi(x)}{2}+\bigl|g(x)-g_{\mathrm{mid}}(x)\bigr|,
	\]
	so that minimizing the mean Hausdorff deviation
	\[
    \int_{X}d_{H}(\{g(x)\},F(x))\,dm(x)
    \]
    amounts to minimizing the $L^{1}$ distance
	from the midpoint. 
    
    The quadratic functional $J$ is preferred because its strict
	convexity yields uniqueness.

	\section{Conclusions}
	\label{sec:final}
	We have studied continuous selections of the interval-valued map
	$F_{A}=[\mu,1-\nu]$ associated with an intuitionistic fuzzy set $A=(\mu,\nu)$, that
	is, ordinary fuzzy membership functions compatible with the intuitionistic
	information. For upper semi-continuous $\mu$ and $\nu$ on a normal space, a
	continuous selection exists by the Kat\v{e}tov--Tong theorem. In the voting model
	of \Cref{sec:voting}, upper semi-continuity is a consequence of the construction,
	and the consensus selection provides an explicit continuous selection, also on
	non-compact universes. On metric spaces, \Cref{thm:lipschitz} characterizes the
	existence of $L$-Lipschitz selections through the non-local inequality
	\eqref{eq:lipschitz-condition} and gives the extremal selections in closed form,
	while \Cref{thm:stability} shows that the set of continuous selections depends in
	a $1$-Lipschitz way on the intuitionistic fuzzy set, with respect to the
	Hausdorff-type distance \eqref{eq:dist-infty}.

	Some questions remain open. Among them, we mention the behaviour of the minimizer
	$g^{*}_{L}$ of \Cref{subsec:criterion} as $L\to\infty$, the existence of a
	distinguished element of the whole set $\mathcal{S}(A)$ of continuous selections,
	which is not compact in general, and the extension of the quantitative results of
	\Cref{sec:lipschitz} to interval-valued intuitionistic fuzzy sets, for which the
	selection problem is genuinely higher-dimensional.

	\section*{Acknowledgements}
    Giuseppe Failla and Marco Elio Tabacchi are supported by funding from the Sustainability Decision Framework
	(SDF) Research Project -- CUP B79J23000540005 -- Grant Assignment Decree No.
	5486 adopted on 2023-08-04. Giuseppe Failla is supported by the Gruppo
	Nazionale per l'Analisi Matematica, la Probabilità e le loro Applicazioni (GNAMPA)
	of the Istituto Nazionale di Alta Matematica (INdAM). Gianmarco La Rosa is supported by the University of Palermo (FFR2026, UNIPA BsD) and also by the Gruppo Nazionale per le Strutture Algebriche, Geometriche
	e le loro Applicazioni (GNSAGA) of the Istituto Nazionale di Alta Matematica (INdAM)
	Francesco Severi.

    \section*{Declarations}

\subsection*{Availability of data and materials}
Not applicable.

\subsection*{Competing interests}
Not applicable. There is no competing interest.

\subsection*{Authors' contributions}
All authors contributed equally to this article.

	\makeatletter
	\renewcommand{\@biblabel}[1]{[#1]\hfill}
	\makeatother
	
\end{document}